\documentclass{article}

\usepackage[utf8]{inputenc}
\usepackage{amsmath,amsfonts,amssymb,amsthm,pb-diagram,picinpar,graphicx,color}
\usepackage{algorithm}
\usepackage{amscd}
\usepackage{algorithmic}
\usepackage{colortbl}
\usepackage{enumerate}
\usepackage{hyperref}
\usepackage{bm}
\usepackage{tikz}
\usepackage{pgf}
\usepackage[T1]{fontenc}
\usepackage{enumitem}
\usepackage{comment}
\usepackage{svg}
\usetikzlibrary{shapes,matrix,calc,arrows,decorations.markings,knots,patterns,intersections,cd}
\usepackage{multirow}
\usepackage{url}

\newcommand{\CC}{\mathbb{C}}
\newcommand{\RR}{\mathbb{R}}
\newcommand{\ZZ}{\mathbb{Z}}
\newcommand{\PP}{\mathbb{P}}
\newcommand{\mcA}{\mathcal{A}}
\newcommand{\mcC}{\mathcal{C}}
\newcommand{\mcB}{\mathcal{B}}

\newcommand{\mcL}{\mathcal{L}}
\newcommand{\mcM}{\mathcal{M}}

\newcommand{\mcR}{\mathcal{R}}
\newcommand{\mcQ}{\mathcal{Q}}
\newcommand{\mcP}{\mathcal{P}}

\newcommand{\bt}{\boldsymbol{t}}
\newcommand{\ba}{\boldsymbol{a}}

\newcommand{\bc}{\boldsymbol{c}}
\newcommand{\bp}{\boldsymbol{p}}
\newcommand{\btau}{\boldsymbol{\tau}}

\newcommand{\MW}{\mathop{\mathrm{MW}}\nolimits}
\newcommand{\Cmb}{\mathop{\mathrm{Cmb}}\nolimits}

\newcommand{\I}{\mathop {\mathrm {I}}\nolimits}

\newcommand{\III}{\mathop {\mathrm {III}}\nolimits}

\newtheorem{thm}{Theorem}[section]
\newtheorem{lem}[thm]{Lemma}     
\newtheorem{cor}[thm]{Corollary}
\newtheorem{prop}[thm]{Proposition}

\theoremstyle{definition}
\newtheorem{defin}[thm]{Definition}

\newtheorem{ex}[thm]{Example}

\theoremstyle{remark}
\newtheorem{rem}[thm]{Remark}

\renewcommand{\thesubparagraph}{\theparagraph.\@arabic\c@subparagraph}
\begin{document}

\begin{center}


{\LARGE \bf 
The realization spaces of certain conic-line arrangements of degree $7$
}

{\large \bf Shinzo Bannai\footnote{Partially supported by JSPS KAKENHI Grant Numbers  JP23K03042},  Hiro-o Tokunaga\footnote{Partially supported by JSPS KAKENHI Grant Number JP20K03561, JP24K06673} and Emiko Yorisaki}


 
\end{center}
\begin{abstract}
We study the embedded topology of certain conic-line arrangements of degree 7.
Two new examples of Zariski pairs are given. Furthermore, we determine the number of connected components of the conic-line arrangements. 
We also calculate the fundamental groups using {\tt SageMath} and the package {\tt Sirocco} in the appendix.
\end{abstract}



\bigskip

{\LARGE \bf Introduction}

\medskip

A collection of a finite number of conics and lines in the
complex projective plane $\PP^2$ is said to be a conic-line 
arrangement (a CL arrangement, for short). If it contains no lines, it is 
said to be a conic arrangement. Compared to line arrangements, which have been
studied by many mathematicians and on which there are
many results from various points of view, there have not been so many
results for CL arrangements up to around 2000 except some results
 on conic arrangements, e.g., \cite{naruki83}.

Since 2000, they have been studied by various mathematicians. For example,
in  
\cite{absst, asstt, amram_garber_teicher07, amram_teicher_uludag03, garber_friedman14, garber_friedman15}, 
the fundamental groups of 
their complements are studied. Also from viewpoints of free divisors, we
find results  such as \cite{schenck2009, dimca_pokora22, 
pokora-szemberg23, macnic}.

In \cite{tokunaga14}, the second author studied the embedded topology of certain
CL arrangements of degree $7$ and gave examples of Zariski pairs for
CL arrangements of degree $7$. He also raised a question 
(see \cite[Remark 6]{tokunaga14}) whether or not 
there exists a Zariski triple for the CL arrangements considered in 
\cite{tokunaga14}. 
In \cite{absst}, another Zariski pair for CL arrangement of degree $7$ was given and
the number of connected components of its realization space was determined.
This article can be considered a continuation of these two articles and we study
the realization spaces of CL arrangements of degree $7$ given by
a similar manner to those in \cite{absst,tokunaga14}.
As a result, we have a negative
answer to the above question in \cite[Remark 6]{tokunaga14} and give
two new Zariski pairs. 
Before explaining our setting and problem explicitly, let us introduce some notation and terminology.

 Let $\mcC\mcL:= \{C_1, \ldots, C_m, L_1, \ldots, L_n\}$ be a CL
 arrangement of $m$-conics and $n$-lines in $\PP^2$.
 By the combinatorics of $\mcC\mcL$ (see \cite{abst2023-1, survey} for 
 the definition of the combinatorics), we mean that of 
the reduced curve $B_{\mcC\mcL}:= \sum_{i=1}^m C_i + \sum_{j=1}^nL_j$ and denote it by $\Cmb_{\mcC\mcL}$.  
 More generally, we denote the combinatorics for a reduced plane curve $B$ by $\Cmb_{B}$.
 The realization space of a given combinatorics $\Cmb_{\mcC\mcL}$ means
 the set of all CL arrangements having the combinatorics $\Cmb_{\mcC\mcL}$
 which we denote by $\mcR(\Cmb_{\mcC\mcL}$). 
 Since all conics and lines are determined by their equations up to 
 non-zero constants, $\mcR(\Cmb_{\mcC\mcL})$ can be regarded as a subset
 of $\PP^{d(d+3)/2}$, where $d = \deg B_{\mcC\mcL}$.
 In this article, we are interested in certain CL arrangements 
$\mcC\mcL$ of degree $7$, which are given in the following way:

\begin{enumerate}
    \item[\rm (i)]  $\mcC\mcL_{ij} = \mcP_i \bigsqcup \mcA_j$ $(i, j = 1, 2)$ where $\mcP_i$ and $\mcA_j$ are subarrangements
    of degree $4$ and $3$ respectively
    such that
      (P1) $\mcP_1 = \{C, L_1, L_2\},\,  \deg C = 2,
    \, \deg L_i = 1 \,
    (i = 1,2)$ with $C \pitchfork (L_1+ L_2)$ and  $C\cap L_1\cap L_2 = \emptyset$,  (P2) $\mcP_2 = \{C_1, C_2\},\,
    \deg C_i = 2 \, (i = 1,2)$ with $C_1\pitchfork C_2$, (A1) $\mcA_1 = \{M_1, M_2, M_3\}$, non-concurrent
    three lines, and (A2) $\mcA_2 = \{D, M\}, \deg D = 2, \deg M =1, D\pitchfork M$.
    We call $\mcP_i$ a plinth for $\mcC\mcL_{ij}$. 

    \item[(ii)] 
    Let $M$ and $D$ be
    a line and a conic in
    $\mcA_j$, respectively. Then any point in $M\cap B_{\mcP_i}$ and $D\cap B_{\mcP_i}$ gives rise to
    a ordinary triple point or
    a tacnode of $M+ B_{\mcP_i}$ and
    $D + B_{\mcP_i}$, respectively.
    
    \item[(iii)] The singularities of 
    $B_{\mcC\mcL_{ij}}$ are at most nodes,
    tacnodes or ordinary triple
    points.

\end{enumerate}

For CL arrangements as above, we have a list as follows:
Here $\Cmb_{ijk}$ denotes the $k$-th combinatorics given by the set
$\mcC\mcL_{ij}$.

\begin{figure}[H]
\centering
\begin{minipage}[b]{0.32\columnwidth}
    \centering
    \includegraphics[width=3.3cm]{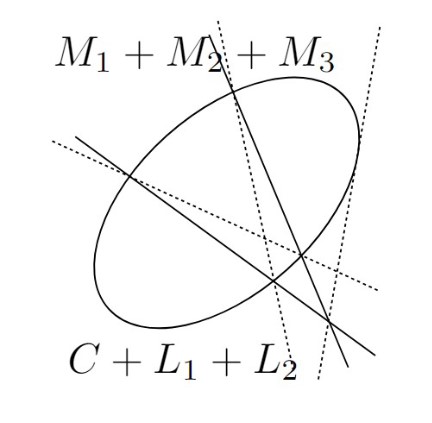}
    \caption{$\Cmb_{111}$}
\end{minipage}
\begin{minipage}[b]{0.32\columnwidth}
    \centering
    \includegraphics[width=3cm]{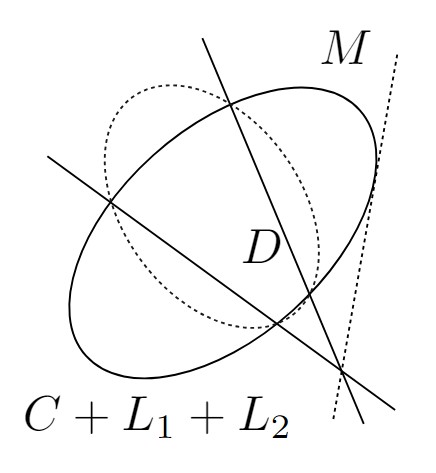}
    \caption{$\Cmb_{121}$}
\end{minipage}
\begin{minipage}[b]{0.32\columnwidth}
    \centering
    \includegraphics[width=3cm]{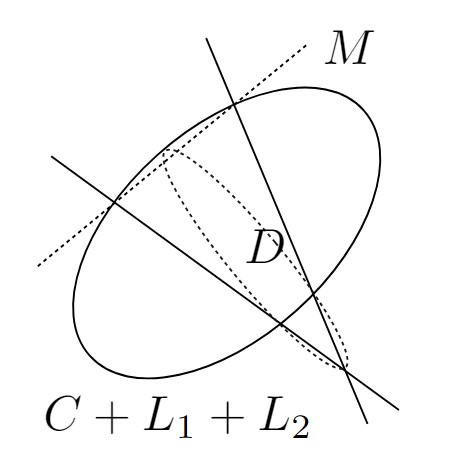}
    \caption{$\Cmb_{122}$}
\end{minipage}
\end{figure}

\begin{figure}[H]
\centering
\begin{minipage}[b]{0.32\columnwidth}
    \centering
    \includegraphics[width=3cm]{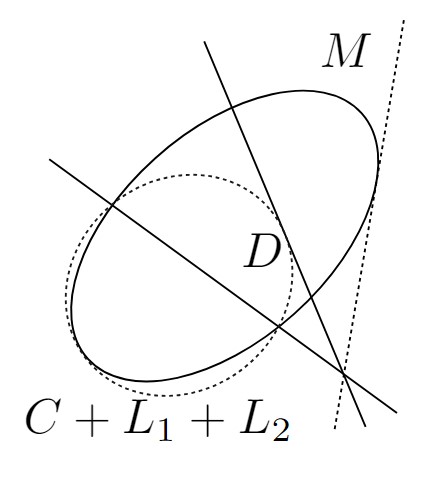}
    \caption{$\Cmb_{123}$}
\end{minipage}
\begin{minipage}[b]{0.32\columnwidth}
    \centering
    \includegraphics[width=3cm]{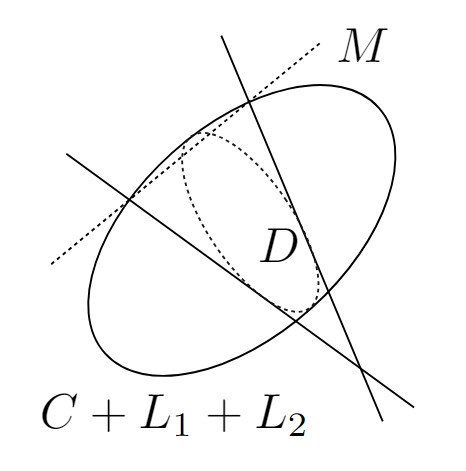}
    \caption{$\Cmb_{124}$}
\end{minipage}
\begin{minipage}[b]{0.32\columnwidth}
    \centering
    \includegraphics[width=3cm]{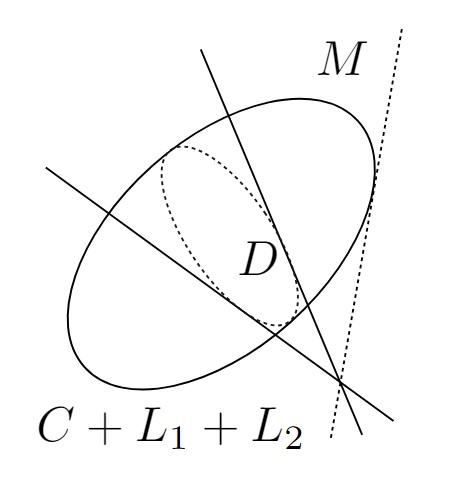}
    \caption{$\Cmb_{125}$}
\end{minipage}
\end{figure}

\begin{figure}[H]
\centering
\begin{minipage}[b]{0.32\columnwidth}
    \centering
    \includegraphics[width=3cm]{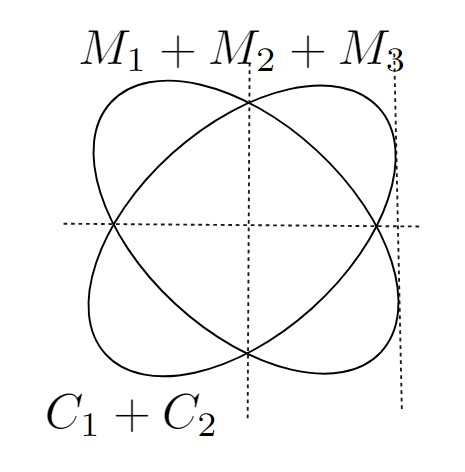}
    \caption{$\Cmb_{211}$}
\end{minipage}
\begin{minipage}[b]{0.32\columnwidth}
    \centering
    \includegraphics[width=3cm]{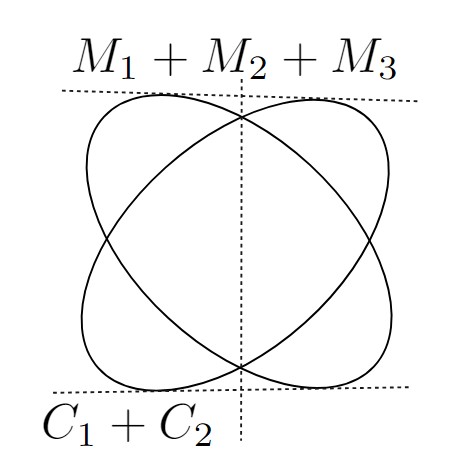}
    \caption{$\Cmb_{212}$}
\end{minipage}
\begin{minipage}[b]{0.32\columnwidth}
    \centering
    \includegraphics[width=3cm]{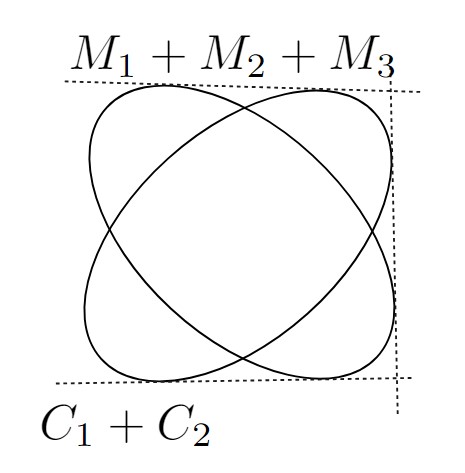}
    \caption{$\Cmb_{213}$}
\end{minipage}
\end{figure}

\begin{figure}[H]
\centering
\begin{minipage}[b]{0.32\columnwidth}
    \centering
    \includegraphics[width=3cm]{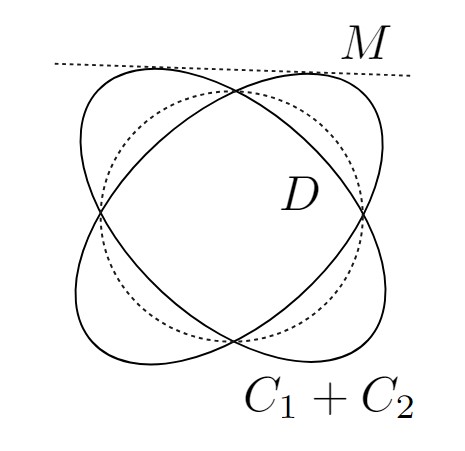}
    \caption{$\Cmb_{221}$}
\end{minipage}
\begin{minipage}[b]{0.32\columnwidth}
    \centering
    \includegraphics[width=3cm]{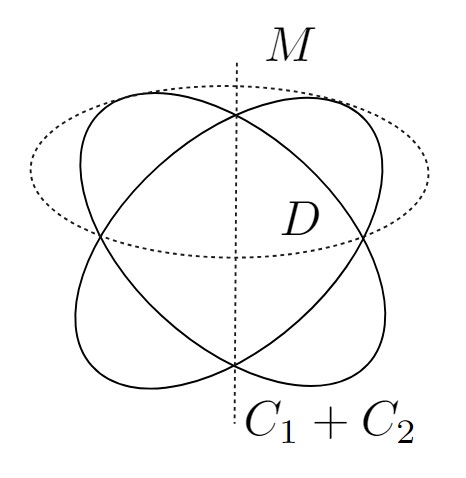}
    \caption{$\Cmb_{222}$}
\end{minipage}
\begin{minipage}[b]{0.32\columnwidth}
    \centering
    \includegraphics[width=3cm]{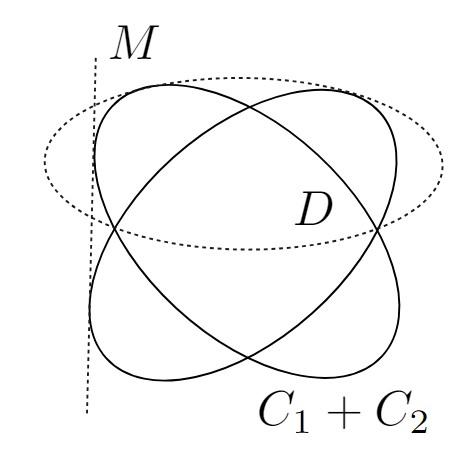}
    \caption{$\Cmb_{223}$}
\end{minipage}
\end{figure}

\begin{figure}[H]
\centering
\begin{minipage}[b]{0.32\columnwidth}
    \centering
    \includegraphics[width=3cm]{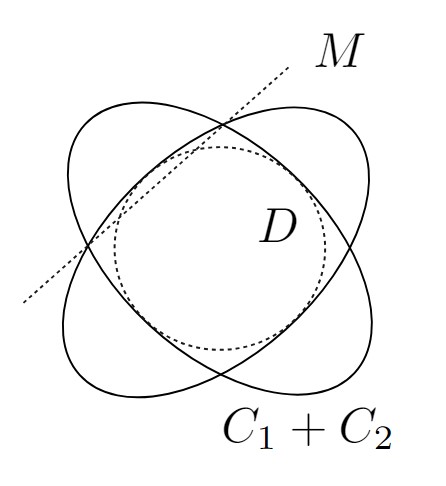}
    \caption{$\Cmb_{224}$}
\end{minipage}
\begin{minipage}[b]{0.32\columnwidth}
    \centering
    \includegraphics[width=3cm]{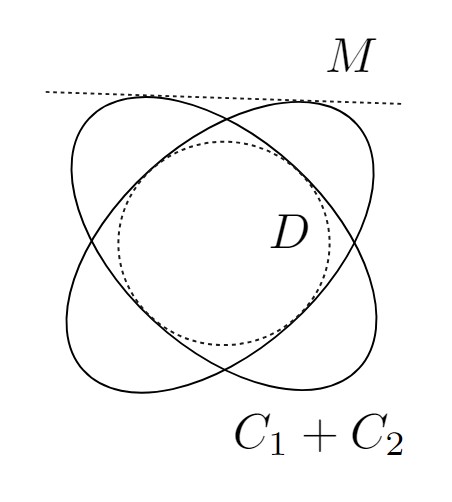}
    \caption{$\Cmb_{225}$}
\end{minipage}
\end{figure}

Fix a conic in $\mcP_i$. As we explain in Subsection \ref{subsec:construction}, some of the above arrangements
are canonically constructed from $5$ points on the conic by using the theory of
rational elliptic surfaces.
 Now our statement is as follows:

 \begin{thm}\label{thm:main}{Let $\Cmb_{ijk}$ denote the combinatorics as in 
 Figures 1-14. Then the following statements hold:
 \begin{enumerate}
     \item[\rm (i)] The space $\mcR(\Cmb_{ijk})$ is connected for $ijk  = 111, 121, 122,125,
     211, 213, 221, 222, 225$.
     \item[\rm (ii)] Each $\mcR(\Cmb_{ijk})$ $ijk = 123, 124, 212, 223, 224$ has two connected components. Moreover, if
     we choose $B_1, B_2 \in \mcR(\Cmb_{ijk})$ so that $B_1$ and $B_2$ belong
     to distinct components, $(B_1, B_2)$ is a Zariski pair.
 \end{enumerate}
 }    
 \end{thm}

\begin{rem}
\begin{itemize}
\item Zariski pairs for the combinatorics $\Cmb_{123}$ and $\Cmb_{212}$ are new, while 
the one for  case $\Cmb_{223}$ was studied in \cite{absst} and the above statements was proved.
\item The case $\Cmb_{224}$ was studied in \cite{tokunaga14} and the existence
of a Zariski triple was expected in \cite[Remark 6]{tokunaga14}.  Theorem~\ref{thm:main}, however,  disproves the existence of such a triple.
\end{itemize}
\end{rem}


\section{Preliminaries}\label{sec:preliminaries}


\subsection{Some rational elliptic surface}\label{subsec:elliptic_surface}
In this article, the theory of elliptic surfaces plays
an important role in both our construction of plane
curves and our proof of Theorem~\ref{thm:main}.
Our main references are \cite{bannai-tokunaga15, BKMT22,  kodaira, miranda-BTES, shioda90, tokunaga14} and we make use of 
results there freely. Here we summarize our convention, notation and
terminology.  For an elliptic surface $\varphi: S \to C$ over a smooth projective curve $C$, we always assume the following:
\begin{enumerate}
    \item[(i)] The fibration $\varphi$ is relatively minimal.
    \item[(ii)] There exists a section $O: C \to S$. We here identify
    $O$ with its image.
    \item[(iii)] There exists at least one singular fiber.
\end{enumerate}

Let $C_o$ be a smooth conic. Choose distinct five points
$z_o, p_1, p_2, p_3$ and $p_4$ on $C_o$. We denote the line
passing through $p_i$ and $p_j$ by $L_{ij}$. Consider a pencil of
conics $\{C_{\lambda}\}_{\lambda \in \Lambda}$ passing through
$p_1, p_2, p_3$ and $p_4$. There exist three distinct values
$\lambda_1, \lambda_2$ and $\lambda_3$ in $\Lambda$ such that
$C_{\lambda_i}$ ($i = 1, 2, 3$) become two distinct lines. We denote
the intersection point between these two lines by $p_0$.
For these values, $C_o + C_{\lambda_i}$ ($i = 1, 2, 3$) 
give rise to conic-line arrangements $\mcP_1$. We may assume that
\[
C_{\lambda_1} = L_{12}+L_{34}, \quad C_{\lambda_2}= L_{13} + L_{24}, \quad
C_{\lambda_3} = L_{14} + L_{23}.
\]
For
other values of $\lambda$, $C_o + C_{\lambda}$ gives rise
to a conic arrangement $\mcP_2$. Put $\mcQ_{\lambda} = C_o + C_{\lambda}, \lambda \in \Lambda$.
Likewise we did in our previous articles \cite{bannai-tokunaga15, BKMT22}, we
associate $(\mcQ_\lambda, z_o)$ with the rational elliptic surface $\varphi_{\mcQ_{\lambda}, z_o} : S_{\mcQ_{\lambda}, z_o} \to \PP^1$, which comes from the double cover $f'_{\mcQ_{\lambda}}: S'_{\mcQ_{\lambda}} \to \PP^2$ branched along $\mcQ_{\lambda}$.
In the following, we always choose $\lambda$ and $z_o$ such that

\begin{center}

$(\ast)$ The tangent line to $C_o$ at $z_o$ meets $C_{\lambda}$ with
two distinct points.

\end{center}

Also the diagram below is the one introduced in \cite{bannai-tokunaga15, BKMT22}
\[
\begin{CD}
S'_{\mcQ_{\lambda}} @<{\mu}<< S_{\mcQ_{\lambda}} @<{\nu_{z_o}}<<S_{\mcQ_{\lambda}, z_o} \\
@V{f'_{\mcQ_{\lambda}}}VV                 @VV{f_{\mcQ_{\lambda}}}V         @VV{f_{\mcQ_{\lambda}, z_o}}V \\
\PP^2@<<{q}< \widehat{\PP^2} @<<{q_{z_o}}< (\widehat{\PP^2})_{z_o},
\end{CD}
\]
where  $\mu$ is the canonical resolution of singularities,  $q$ is a composition of a finite number of blowing-ups so that the branch locus becomes smooth 
and $f_{\mcQ_{\lambda}}$ is the induced double cover. 
The pencil of lines through $z_o$ gives rise a pencil $\Lambda_{z_o}$ of curves of genus $1$. We denote the resolution of indeterminacy
of $\Lambda_{z_o}$ by $\nu_{z_o}$ and $q_{z_o}$ is a composition of two blowing-ups induced by $\nu_{z_o}$.  We also have an induced double cover $f_{\mcQ_{\lambda}, z_o} : S_{\mcQ_{\lambda}, z_o} \to (\widehat{\PP^2})_{z_o}$. 
The generic fiber $E_{\mcQ_{\lambda},z_o}$ can be consider an elliptic curve over $\CC(\PP^1)(\simeq \CC(t))$.
The induced double cover  $f_{\mcQ_{\lambda}, z_o}$ 
coincides  with the quotient morphism determined by the involution $[-1]$ on $S_{\mcQ_{\lambda},z_o}$, which is given by the inversion with respect to the group law on $E_{\mcQ_{\lambda},z_o}$.  
Let $E_{\mcQ_{\lambda},z_o}(\CC(t))$ be the set of $\CC(t)$-rational points of $E_{\mcQ_{\lambda},z_o}$
and let $\MW(S_{\mcQ_{\lambda},z_o})$ be the set of sections. 
By an
integral section, we mean a suction $s$ with $s\cdot O = 0$.
 In \cite{shioda90}, Shioda defined a ${\mathbb{Q}}$-valued
bilinear form $\langle\,\, , \,\,\rangle$  on  $E_{\mcQ_{\lambda},z_o}(\CC(t))$ called the height pairing, by which the free part of $E_{\mcQ_{\lambda},z_o}(\CC(t))$ becomes a lattice. We make use of this lattice structure in order to find
elements in $\mcA_j$ ($j = 1,2$). When we describe $E_{\mcQ, z_o}(\CC(t))$,
we take this structure into account.
It is known that there is a bijection between $\MW(S_{\mcQ_{\lambda},z_o})$ and $E_{\mcQ_{\lambda},z_o}(\CC(t))$.
For $s \in \MW(S_{\mcQ_{\lambda},z_o})$, we denote the rational point corresponding to $s$ by $P_s$, and for $P \in E_{\mcQ_{\lambda},z_o}(\CC(t))$, we denote the section corresponding to $P$ by $s_P$.
Under this correspondence, we have $s_{P_1\dot{+}P_2} = s_{P_1}\dot{+}s_{P_2}$.
We also write $\mcC_s:= f'_{\mcQ_{\lambda}} \circ \mu \circ \nu_{z_o}(s) \subset \PP^2$ for a section $s \in \MW(S_{\mcQ_{\lambda},z_o})$.


Here are some properties of $\varphi_{\mcQ_{\lambda}, z_o}: S_{\mcQ_{\lambda}, z_o} \to \PP^1$ (See \cite{bannai-tokunaga15, BKMT22, tokunaga14, oguiso-shioda}):

\underline{The Case $\lambda = \lambda_1, \lambda_2, \lambda_3$}

\begin{itemize}
    \item  There exist $6$ singular fibers for $\varphi_{\mcQ_{\lambda}, z_o}$.
    All of them are of types $\I_2$. 
    They arise from
    the tangent line $l_{z_o}$ at $z_o$ and lines connecting $z_o$ and $p_i$
    ($0 \le i \le 4$). We denote them by $F_{\infty}$ and $F_i$ ($0 \le i \le 4$), respectively,  and their irreducible decomposition by $F_{\bullet} = 
    \Theta_{\bullet, 0} + \Theta_{\bullet,1}$ $\bullet = \infty, 0, 1, \ldots, 4$.
    \item The group $E_{S_{\mcQ_{\lambda}, z_o}}(\CC(t))$ is isomorphic
    to $(A_1^{\ast})^{\oplus 2} \oplus (\ZZ/2\ZZ)^{\oplus 2}$.

    \item In order to describe explicit generators of $E_{\mcQ_{\lambda}, z_o}(\CC(t))$, we consider the case $\lambda = \lambda_1$. In this case, 
    $C_o$ and $L_{ij}$ ($1\le i < j \le 4$) give rise to elements of
    $E_{\mcQ_{\lambda}, z_o}(\CC(t))$ as follows:
    \begin{enumerate}
        \item[(i)] $C_o$, $L_{12}$ and $L_{34}$ give rise $2$-torsions, 
        which we denote by
        $P_{C_o}$, $P_{12}$ and $P_{34}$, respectively. Note that
        $P_{C_o} = [-1]P_{C_o}$, $P_{12} = [-1]P_{12}$ and $P_{34} = [-1]P_{34}$.
        \item[(ii)] For each $(i, j) \in \{(1, 3), (1,4), (2, 3), (2, 4)\}$,
        $L_{ij}$ gives rise to two points $[\pm]P_{ij} \in E_{\mcQ_{\lambda},z_o}(\CC(t))$.
        \item[(iii)] We may assume that the free part of $E_{\mcQ_{\lambda},z_o}(\CC(t))$ generated by $P_{13}$ and $P_{14}$, i.e.,
        \[
        (A_1^{\ast})^{\oplus 2} \cong \ZZ P_{13}\oplus \ZZ P_{14}
        \]
        and $P_{23} = P_{14}\dot{+}P_{C_o}, P_{24} = P_{13} \dot{+}P_{C_o}$.
        \item[(iv)] For each $(i, j) \in \{(1, 3), (1,4), (2, 3), (2, 4)\}$,
        $\mcC_{[2]P_{ij}} =  \mcC_{[-2]P_{ij}}$ is a conic inscribed by 
        $\mcQ_{\lambda_1}$ such that $z_o \in \mcQ_{\lambda_1}\cap \mcC_{[2]P_{ij}}$.
    \end{enumerate}
    
\end{itemize}

\underline{The Case $\lambda \neq \lambda_1, \lambda_2, \lambda_3$}

\begin{itemize}
\item There exist $5$ reducible singular fibers. 
 All of them are of types either $\I_2$ or $\III$. They arise from
    the tangent line $l_{z_o}$ at $z_o$ and lines through $z_o$ and $p_i$
    ($1 \le i \le 4$). We denote them by $F_{\infty}$ and $F_i$ ($1 \le i \le 4$), respectively,  and their irreducible decomposition by $F_{\bullet} = 
    \Theta_{\bullet, 0} + \Theta_{\bullet,1}$ $\bullet = \infty, 1, \ldots, 4$.

\item The group $E_{S_{\mcQ, z_o}}(\CC(t))$ is isomorphic to  $(A_1^*)^{\oplus 3}\oplus \ZZ/2\ZZ$. 
The unique $2$-torsion point arises from $C_o$, which we denote by $P_{C_o}$. 

\item Each $L_{ij}$ gives two elements in $E_{\mcQ_{\lambda}, z_o}$ and we
denote them by  $[\pm 1]P_{ij}$, which satisfy the following properties:
\begin{enumerate}
  \item[(i)] Since $\langle P_{1j}, P_{1j} \rangle = 1/2$ ($2\le j \le 4$), $
\langle P_{1j}, P_{1k}\rangle = 0$ ($2\le j < k \le 4$), 
we may assume
\[
(A_1^*)^{\oplus 3} \cong \ZZ P_{12} \oplus \ZZ P_{13} \oplus \ZZ P_{23}
\]
and  $P_{ij} \dot {+} T = P_{kl}$, where $\{i, j, k, l\} =\{ 1, 2, 3, 4\}$.

 \item[(ii)] $\mcC_{[2]P_{ij}} =  \mcC_{[-2]P_{ij}}$ is a conic inscribed by 
        $\mcQ_{\lambda}$ such that $z_o \in \mcQ_{\lambda}\cap \mcC_{[2]P_{ij}}$.
\end{enumerate}        
\end{itemize}

\subsection{Construction of  lines and conics in \texorpdfstring{$\mcA_j\, (j = 1,2)$}{} via 
\texorpdfstring{$S_{\mcQ_\lambda, z_o}$}{}}\label{subsec:construction}

We here explain our method in constructing lines and conics in 
$\mcA_j$ ($j = 1,2$). This method plays a crucial role to 
consider a member of $\mcR(\Cmb_{ijk})$.
Choose $P \in E_{\mcQ_{\lambda}, z_o}$ and let $s_P$ be the section.
In \cite{masuya24}, Masuya introduced  a  {\it line point}  as follows:

\begin{defin}\label{def:line_point}{\rm 
$P$ is said to be a line-point if $\tilde{f}_{\mcQ_{\lambda}, z_o}(s_P)$
is a line. Also a section $s \in \MW(S_{\mcQ_\lambda, z_o})$ is said to be
a line-section if $\tilde{f}_{\mcQ_\lambda, z_o}(s)$ is a line.
}    
\end{defin}

Any line-point is characterized by the following lemma:

\begin{lem}\label{lem:line_point}{{\rm (\cite[Lemma 9]{bannai-tokunaga17})} Let
$s \in \MW(S_{\mcQ_{\lambda}, z_o})$ be an integral section  with $s\cdot \Theta_{\infty,1}
= 1$. Then $\tilde{f}_{\mcQ_{\lambda}, z_o}(s)$ is a line $L_s$ such that
\begin{enumerate}
    \item[\rm (i)] the intersection multiplicity at every intersection point between
    $L_s$ and $\mcQ_{\lambda}$ is even,
    \item[\rm (ii)] $z_o \not\in L_s$.
\end{enumerate}
Conversely, any line satisfying the above two conditions gives rise to two
sections $s_{L^{\pm}}$ such that $s_{L^{\pm}}\cdot O = 0$ and
$s_{L^{\pm}}\cdot \Theta_{\infty, 1}=1$.
}
\end{lem}

As for an integral section $s$ with $s\cdot \Theta_{\infty, 0} = 1$, we have the 
following lemma:

\begin{lem}\label{lem:conic_point}{{\rm \cite[Lemma 2.12]{masuya24}} Let $s \in 
\MW(S_{\mcQ_{\lambda}, z_o})$ be an integral section with $s\cdot\Theta_{\infty, 0}$.
Then $\tilde{f}_{\mcQ_{\lambda}, z_o}(s)$ is a smooth conic satisfying either
\begin{enumerate}
    \item[\rm (i)] $\tilde{f}_{\mcQ_{\lambda}, z_o}(s)$ is the irreducible component
     of $\mcQ_{\lambda}$ through $z_o$, or 
    \item[\rm (ii)] $\tilde{f}_{\mcQ_{\lambda}, z_o}(s)$ is tangent to $\mcQ_{\lambda}$ at
    $z_o$ and the intersection multiplicity at every intersection point between
    $\tilde{f}_{\mcQ_{\lambda}, z_o}(s)$ and $\mcQ_{\lambda}$ is even.
\end{enumerate}
Conversely, any conic C that satisfies one of the above conditions gives rise
to two integral sections $s_{C^{\pm}}$ such that $s_{C^{\pm}}\cdot \Theta_{\infty, 0} = 1$.
}
\end{lem}
Note that although the "converse" part of Lemma~\ref{lem:conic_point} was not given in 
\cite{masuya24}, it follows from our construction of $S_{\mcQ_{\lambda}, z_o}$.
By Lemmas~\ref{lem:line_point} and \ref{lem:conic_point}, we see that lines and conics
in $\mcA_j$ $(j = 1, 2)$ are canonically obtained from sections described as above and 
vice versa. Let $s$ be a section in $\MW(S_{\mcQ_{\lambda}, z_o})$ and let $P_s$ be the
corresponding point in $E_{\mcQ_{\lambda}, z_o}$. If we choose a section as
in Lemmas~\ref{lem:line_point} and \ref{lem:conic_point}, the we have the table below:

\begin{table}[htbp]
\centering
\begin{tabular}{|c|c|}\hline
Line or conic & Points in $E_{\mcQ_{\lambda}, z_o}(\CC(t))$  \\ \hline
$L_0$ &
    $P_{13} \dot{\pm} P_{14} \dot{+} P_{34}$
    \\ \hline
$D(1,1)$ &
    $P_{13} \dot{\pm} P_{12}$
    , $P_{13} \dot{\pm} P_{34}$
    , $P_{14} \dot{\pm} P_{12}$
    , $P_{14} \dot{\pm} P_{34}$
    \\ \hline
$D(1,2)$ &
    $P_{13} \dot{\pm} P_{14}$
    , $P_{13}\dot{\pm}P_{23}$
    \\ \hline
$D(1,4)$ &
    $[2]P_{13}$
    , $[2]P_{14}$
    \\ \hline
$L_{b}$ &
    $P_{12} \dot{\pm} P_{13} \dot{\pm} P_{23}$
    \\ \hline
$D(2,2)$ &
    $P_{12} \dot{\pm} P_{23}$
    , $P_{12} \dot{\pm} P_{13}$
    , $P_{13} \dot{\pm} P_{23}$
    , $P_{12}\dot{\pm}P_{14}$
    , $P_{13}\dot{\pm}P_{14}$
    , $P_{12}\dot{\pm}P_{24}$
    \\ \hline
$D(2,4)$ &
    $[2]P_{12}$
    , $[2]P_{13}$
    , $[2]P_{23}$
    \\ \hline
\end{tabular}
\caption{Line points, conic points and their corresponding curves}
\label{tab:generator}
\end{table}
\begin{rem}{\rm \quad 
\begin{itemize}
 \item  As $\tilde{f}_{\mcQ_{\lambda}, z_o}(s_P) 
= \tilde{f}_{\mcQ_{\lambda}, z_o}(s_{[-1]P})$, we only give one of the corresponding
two points. Also, $\dot{\pm}$ can be chosen freely.
  \item Since $P_{ij} = P_{kl} \dot{+}T$, relations $P_{12}\dot{\pm}P_{23} = P_{34}\dot{\pm}P_{14}$,
  $P_{12}\dot{\pm}P_{14} = P_{34}\dot{\pm}P_{23}$, etc hold. 
 \end{itemize}   
 }
\end{rem}


\section{Approach to construct  CL arrangements with prescribed \texorpdfstring{$\Cmb_{ijk}$}{}}\label{subsec:explicit-construction}

In this section, we give rough ideas for explicit construction
of plane curves with prescribed
$\Cmb_{ijk}$.
We first give a combinatorial classification of lines and a conic
in $\mcA_j$ ($j = 1,2$).


\underline{The case of $\mcP_1$:} We may assume
$L_{12}, L_{34} \in \mcP_1$
and put $L_1 = L_{12}, L_2= L_{34}$.

Let $M$ be a line
in $\mcA_j$ ($j = 1, 2)$. By our assumption (ii) in the
Introduction, we infer that $M$ is $L_{13}, L_{14}, L_{23}, L_{24}$ or line $L_0, L'_0$
through $p_0$ and tangent to $C$.

Let $D$ be the smooth conic in $\mcA_2$. Again by our assumption (ii), we infer that
$D$ is one of conics of type $D(1,j)$ ($j = 0, 1, 2, 4$) as follows:
 \begin{enumerate}
    \item[(a)] $D(1,0)$ passes through $p_1, p_2, p_3$ and $p_4$.
    \item[(b)] $D(1,1)$ passes through $p_0, p_i, p_j$ ($i \in \{1,2\}, j\in \{3, 4\}$) and
    tangent to $C$.
    \item[(c)] $D(1,2)$ passes through $p_1$ and $p_2$ (resp. $p_3$ and $p_4$) and
    tangent to $L_{34}$ (resp. $L_{12}$) and $C$.
    \item[(d)] $D(1,4)$ is a conic inscribed by $B_{\mcP_1}$.
 \end{enumerate}

\underline{The case of $\mcP_2$:} Let $M$ be a line in $\mcA_j$ ($j=1, 2$). By
assumption (ii) in the Introduction, we infer that $M$ is a bitangent line $L_{b}$
to $B_{\mcP_2}$ or $L_{ij}$ ($1 \le i < j \le 4$). Note that there exist four bitangent
lines to $B_{\mcP_2}$.

Let $D$ be the conic in $\mcA_2$. By assumption (ii), we infer that
$D$ is one of conics of type $D(2,j)$ ($j= 0, 2, 4$) as follows:
 \begin{enumerate}
     \item[(a)] $D(2,0)$ passes through $p_1, p_2, p_3$ and $p_4$.
     \item[(b)] $D(2,2)$ passes through $p_i, p_j$ and tangent to both $C_1$ and $C_2$.
     \item[(c)] $D(2,4)$ is tangent to $B_{\mcP_2}$ at $4$ distinct points.
 \end{enumerate}

We next explain our setting about an explicit Weierstrass equation for $E_{\mcQ_{\lambda}, z_o}$, which
gives an affine equation of $S_{\mcQ_{\lambda}, z_o}$.
This setting plays an important role
to prove Theorem~\ref{thm:main}.
Let $C_o$ and $z_o$, $p_1, \ldots, p_4$ be the smooth conic and 
distinct $5$ points on it as in Subsection~\ref{subsec:elliptic_surface}.
We take
homogeneous coordinate $[T, X, Z]$ of $\PP^2$ such that $C_o$ is given by
$XZ - T^2 = 0$ and $z_o=[0,1,0]$. Note that $l_{z_o}$ is given by $Z=0$. Let $(t, x)$, $t = T/Z, x = X/Z$ be affine coordinates
of $\PP^2\setminus l_{z_o}$. 
Put $p_i:=(t_i,t_i^2)$, $t_i \in \mathbb{C}$  $(i=1,2,3,4)$.
Let $\bt=(t_1,t_2,t_3,t_4)$ and $\lambda \in \CC$.
We define $\mcM$ as follows:
\[
\mcM:= \{\btau= (\lambda, \bt) \in \CC \times \CC^4 \mid t_i \neq t_j \,\, (i \neq j) \}.
\]
Under these settings, let  $C_{\btau}$  be the conic given by 
\[
C_{\btau}: c_{\btau}(t, x):= \lambda(x - t^2)+(x - (t_1+ t_2)t + t_1t_2)(x - (t_3 + t_4)t + t_3t_4),
\quad (\lambda,\bt) \in \mcM.
\]
With this equation, $E_{\mcQ_{\btau}, z_o}$ is given by
the Weierstrass equation:
\[
E_{\mcQ_{\btau}, z_o}: y^2 = f_{\btau}(t, x), \quad
f_{\btau}(t, x) = (x - t^2)c_{\btau}(t,x)
\]

Note that $\mcQ_{\btau}$ and explicit generators of 
$E_{\mcQ_{\btau}, z_o}(\CC(t))$ are
determined by $C_o$ and $z_o, p_1, \ldots, p_4$ by Subsection~\ref{subsec:elliptic_surface}.

Now let us explain how we construct CL arrangement with $\Cmb_{ijk}$. First we may assume that quartics $B_{\mcP_1}, B_{\mcP_2}$ are given
by $\mcQ_{\btau}$ defined by the equation of
the form $
f_{\btau}(t, x) = 0,\, \btau \in \mcM$.
 Also we keep our notation for lines and conics in $\mcA_j$ ($j = 1, 2$) as in the beginning of this section.

The quartics  $B_{\mcP_1}$ are given by three values $\lambda_1 = 0, \lambda_2=-(t_1 - t_4)(t_2 - t_3)$
and $\lambda_3 = -(t_1-t_3)(t_2- t_4)$ for a  fixed $\bt$ and we have $\mcP_1 = \{C_o, L_{12}, L_{34}\},
 \{C_o, L_{13}, L_{24}\}$ and $ \{C_o, L_{14}, L_{23}\}$, respectively. 
On the other hand, once we choose one of three values
$\lambda_1, \lambda_2, \lambda_3$ and fix it, by interchanging the coordinates
 of $\bt$ continuously,
pair of lines $\{L_{ij}, L_{kl}\}$ $i< j, k<l, \{i, j,k, l\} = \{1, 2, 3, 4\}$ are
also continuously interchanged.
Hence we may assume that $\mcP_1$ is given by $\lambda_1 = 0$ and 
$B_{\mcP_1} = \{C_o, L_{12}, L_{34}\}$.
Here are some more remarks:

\begin{itemize}
 \item $B_{\mcP_1}$ is determined by a $2$-partition of $\{1,2,3,4\}$.
 Since $p_0$ is determined by this $2$-partition, two tangent lines to $C_o$ that pass through
 $p_0$ is also canonically determined by $\bt$.
 \item Fix $\bt$. Then, smooth conics of type $D(1,0), D(2,0)$ which passes through $p_1, p_2, p_3$ and $p_4$ are give by $C_{\lambda}$
 for some $\lambda$.
 \item  In order to describe $\mcR(\Cmb_{ijk})$, we define two disjoint subsets, $\mcM_1$ and $\mcM_2$, of $\mcM$ as follows:
\[
\mcM_1:= \{\btau = (\lambda_1,\bt) = (0, \bt) \in \mcM\},  
\quad 
\mcM_2:= \{\btau= (\lambda, \bt) \in \mcM \mid \lambda \neq \lambda_1, \lambda_2, \lambda_3\}.
\]
\end{itemize}

In the following, we explain how we construct conic-line arrangements with $\Cmb_{ijk}$.

For $\mcP_1$, we may assume that $\mcP_1 = \{C_o, L_{12}, L_{34}\}$ for some fixed $\btau=(0,\bt) \in \mcM_1$.

\underline{$\Cmb_{111}$}:  In this case, $\mcA_1$ is one of the following
\[
\{L_{13}, L_{24}, L_0\}, \{L_{13}, L_{24}, L'_0\}, \{L_{14}, L_{23}, L_0\}, \{L_{14}, L_{23}, L'_0\}
\]

\underline{$\Cmb_{121}$}: Any conic of type $D(1,0)$ is given by $c_{\btau'}(t, x) = 0$ for some $\btau'=(\lambda',\bt) \in \mcM_2$, which we denote by $D$.
Hence we may assume $\mcA_2$ is given by $\{L_0, D\}$ or $\{L'_0, D\}$.

\underline{$\Cmb_{122}$}: The conic in $\mcA_2$ is of type $D(1,1)$.
We denote a conic of type $D(1,1)$ passing through $p_i, p_j$ by $D_{ij}$.
Hence $\mcA_2$ is one of the following:
\[
\{L_{13}, D_{24}\}, \{L_{14}, D_{23}\}, \{L_{23}, D_{14}\}, \{L_{24}, D_{13}\}.
\]
We may assume that $D_{ij}$ is tangent to $C_o$ at $z_o$.
Then $D_{ij}$ is a parabola through $p_0, p_i$ and $p_j$, which is uniquely determined. Hence any CL
arrangements with $\Cmb_{122}$ is determined by $p_1, p_2, p_3, p_4$ and $z_o$.

\underline{$\Cmb_{123}$}: 
The conic in $\mcA_2$ is of type $D(1,2)$.
We denote a conic of type $D(1,2)$ passing through $p_i, p_j$ by $D_{ij}$.
We may assume that $D_{ij}$ is tangent to $C_o$ at $z_o$. Hence we see that $\mcA_2$ is one of the following:
\[
\{L_0, D_{12}\},
\{L'_0, D_{12}\}, \{L_0, D_{34}\},
\{L'_0, D_{34}\}.
\]
Note that $L_0, L'_0, D_{ij}$ is obtained from sections of
$S_{\mcQ_{\btau}, z_o}$ as in the table in Subsection~\ref{subsec:construction}.
Hence, every CL arrangement with $\Cmb_{123}$ is 
is determined by $p_1, p_2, p_3, p_4$ and $z_o$.

\underline{$\Cmb_{124}$}
The conic in $\mcA_2$ is of type $D(1,4)$, which we denote by $D$.
We may assume that $D$ is tangent to $C_o$ at $z_o$.
Hence we see that $\mcA_2$ is one of the following:
\[
\{L_{13}, D\},
\{L_{14}, D\}, \{L_{23}, D\},
\{L_{24}, D\}.
\]
Note that $L_{ij}, D$ as above are obtained from sections of
$S_{\mcQ_{\btau}, z_o}$ as in the table in Subsection~\ref{subsec:construction}.
Hence, every CL arrangement with $\Cmb_{124}$ is 
is determined by $p_1, p_2, p_3, p_4$ and $z_o$.

\underline{$\Cmb_{125}$}: The conic in $\mcA_2$ is of type $D(1,4)$, which we denote by $D$.
We may assume that $D$ is tangent to $C_o$ at $z_o$.
Hence we see that $\mcA_2 = \{L_0, D\}, \{L'_0, D\}$.
Note that $L_0, L'_0, D$ is obtained from sections of
$S_{\mcQ_{\btau}, z_o}$ as in the table in Subsection~\ref{subsec:construction}.
Hence, every CL arrangement with $\Cmb_{125}$ is determined by $p_1, p_2, p_3, p_4$ and $z_o$.

For $\mcP_2$, we may assume that $\mcP_2 = \{C_o, C_{\btau} \}$ for some fixed $\btau \in \mcM_2$.
There exist four bitangent lines for $B_{\mcP_2}$, which we denote by $L_{b_i}$ ($1\le i \le 4$). Table 1 in the previous
subsection shows that all
bitangent lines of $B_{\mcP_2}$ are determined by  $p_1, p_2, p_3, p_4$ and are
canonically constructed if we choose $z_o$.

\underline{$\Cmb_{211}$}: $\mcA_1$ consists of three lines as follows:
$\{L_{12}, L_{34}, L_b\}, \{L_{13}, L_{24}, L_b\}, \{L_{14}, L_{23}, L_b\}$, where
there are $4$ possibilities for $L_b$. By Table 1, every bitangent lines are given by three line points. Hence, every CL arrangement with $\Cmb_{211}$ is
determined by  $\btau$. 

\underline{$\Cmb_{212}$}: $\mcA_1$ consists of pair of four bitangent lines and $L_{ij}$.
There are $24$ possibilities for such collections. 
Yet, likewise $\Cmb_{211}$, every CL arrangement with $\Cmb_{212}$ is determined by $\btau$. 

\underline{$\Cmb_{213}$}: $\mcA_1$ consists of three of four bitangent lines.
Likewise $\Cmb_{211}$, every CL arrangement with $\Cmb_{213}$ is
determined by $\btau$. 

\underline{$\Cmb_{221}$}: $\mcA_2$ consists of a bitangent line and 
a smooth conic of type $D(2,0)$.
Every CL arrangement with $\Cmb_{221}$ is determined by $\btau$ and 
another value $\lambda'$. 

\underline{$\Cmb_{222}$}: $\mcA_2$ consists of a line $L_{ij}$ and a conic of type $D(2,2)$.
We denote a conic of type $D(2,2)$ passing through $p_i, p_j$ by $D_{ij}$. 
Then we infer that $\mcA_2$ is of the form $\{L_{ij}, D_{kl}\}$, $\{i,j,k,l\} = \{1,2,3,4\}$.
We may assume that $D_{kl}$ is tangent to $C_o$ at $z_o$.
By Table 1, every CL arrangement with $\Cmb_{222}$ is constructed from  $\btau$ 
and $z_o$
by choosing $P_{ij}$'s appropriately.

\underline{$\Cmb_{223}$}:
$\mcA_2$ consists of a bitangent line $L_{b}$ and a conic of type $D(2,2)$.
We denote a conic of type $D(2,2)$ passing through $p_i, p_j$ by $D_{ij}$. 
Then we infer that $\mcA_2$ is of the form $\{L_{b}, D_{ij}\}$.
We may assume that $D_{ij}$ is tangent to $C_o$ at $z_o$.
By Table 1, every CL arrangement with $\Cmb_{223}$ is constructed from $\btau$ 
and $z_o$ by choosing $P_{ij}$'s appropriately.

\underline{$\Cmb_{224}$}:
$\mcA_2$ consists of a line $L_{ij}$ and a conic $D$ of type $D(2,4)$, which we denote by $D$.
We may assume that $D$ is tangent to $C_o$ at $z_o$.
By Table 1, every CL arrangement with $\Cmb_{224}$ is constructed from $\btau$ 
and $z_o$ by choosing
$P_{ij}$'s appropriately.

\underline{$\Cmb_{225}$}: $\mcA_2$ consists of a bitangent line $L_b$ and a conic of type $D(2,4)$, which we denote by $D$.
We may assume that $D$ is tangent to $C_o$ at $z_o$. 
By Table 1, every CL arrangement with $\Cmb_{225}$ is constructed from $\btau$ 
and $z_o$ by choosing $P_{ij}$'s appropriately.

\section{Proof of Theorem~\ref{thm:main}}

\subsection{Our strategy}\label{subsec:strategy}

Let us explain our strategy to prove Theorem~\ref{thm:main}.  Our approach is similar to that we take in 
\cite[Section 3]{absst}. As we see in Section~\ref{subsec:explicit-construction}, every
CL arrangement with $\mcP_i$ is obtained from $\btau$, $z_o$
and another parameter $\lambda'$. Once we fixed, we have construct CL arrangement with $\Cmb_{ijk}$  in a canonical way via
sections of $S_{\mcQ_{\lambda}, z_o}$ except those involving conics of type $D(1,0)$ and $D(2,0)$. 
Our basic idea to prove Theorem~\ref{thm:main} is to  CL arrangements with fixed $\Cmb_{ijk}$ by moving $\btau$,
which is done in \cite[Lemma 3.1, Remark 3.2, Corollary 3.3]{absst}.
Let us explain it more precisely. Let $\mcM$, $\mcM_1$ and $\mcM_2$ be those in Section~\ref{subsec:explicit-construction}.
We first choose $z_o$.
Since $p_i = (t_i, t_i^2)$, to move points $(p_1, p_2, p_3, p_4)$ and $\lambda$ can be considered as to move $\btau = (\lambda,\bt)$ in $\mcM$ along a path in $\mcM$ as in Section~\ref{subsec:explicit-construction}.
We consider $\gamma : [0,1] \to \mcM$, $s \mapsto \gamma(s) = (\lambda(s),(t_1(s),t_2(s),t_3(s),t_4(s)))$ as such a path.
Let $\btau_o=(\lambda_o,(a_1,a_2,a_3,a_4)),\btau_o'=(\lambda_o',(a_1',a_2',a_3',a_4')) \in \mcM$, $\gamma(0)=\btau_o$ and $\gamma(1)=\btau_o'$.
If $t_i(0)=a_i$ and $t_i(1)=a_j'$, we say \lq\lq $a_i$ goes to $a_j'$ along a path in $\mcM$ \rq\rq  and denote it by $a_i \rightsquigarrow a_j'$.
In our proof of Theorem~\ref{thm:main}, we exploit $\mcM_1$ to describe connected components of $\mcR(\Cmb_{ijk})$ for $ijk = 111, 122, 123, 124, 125$,
while we exploit $\mcM_2$ to describe connected components of  $\mcR(\Cmb_{ijk}))$ $ijk = 121, 211, 212, 213, 224$.
In this section, we keep the same notation for lines and conics as 
those given in Section \ref{subsec:explicit-construction}.
Now we prove Theorem~\ref{thm:main} based on case-by-case arguments.

\subsection{\texorpdfstring{$\Cmb_{111}$}{}}\label{subsec:111}
Let $C_o$ be the conic as before and choose $\btau=(\lambda_1, \bt) \in \mcM_1$.
We put two elements $B_{\btau}, B'_{\btau} \in \mcR(\Cmb_{111})$ as follows:
\begin{eqnarray*}
   \mcQ_{\btau} & = & C_o + L_{12, \btau} + L_{34, \btau} \\
   B_{\btau}&:=& \mcQ_{\btau} + L_{i_1j_1, \btau} + L_{i_2j_2, \btau} + L_{0, \btau}, \\
B'_{\btau}&:= & \mcQ_{\btau}+ L_{i_1j_1, \btau} + L_{i_2j_2, \btau} + L'_{0, \btau}. 
\end{eqnarray*}
Note that once we choose $L_{i_1j_1, \btau}$, $L_{i_2j_2, \btau}$ is automatically determined.

Choose $\ba = (-2,-1, 1, 2)$, $\btau_o=(\lambda_1,\ba) \in \mcM$.
In this case, $L_{12, \btau_o}: x + 3t + 2 = 0$, $L_{34, \btau_o}: x - 3t+2=0$,
$L_{13, \btau_o}: x + t -2 = 0$, $L_{24, \btau_o}: x - t -2 = 0$ and $p_0 = (0, -2)$. Hence we put
\[
L_{0, \btau_o}: x +2\sqrt{2}t +2 = 0, \quad L'_{0,\btau_o}: x -2\sqrt{2}t +2 = 0.
\]
Define $B_{\btau_o}, B'_{\btau_o} \in \Cmb_{111}$ to be
\begin{eqnarray*}
  B_{\btau_o} &:= & \mcQ_{\btau_o}+ L_{13, \btau_o} + L_{24, \btau_o} + L_{0,\btau_o}, \\
  B'_{\btau_o} &:= &  \mcQ_{\btau_o} + L_{13, \btau_o} + L_{24, \btau_o} + L'_{0,\btau_o}.
\end{eqnarray*}
Note that $B'_{\btau_o}$ is transformed to $B_{\btau_o}$ by $(t, x) \mapsto (-t, x)$.
Now choose $B\in \mcR(\Cmb_{111})$ arbitrary. By taking suitable coordinates of $\PP^2$ so that the conic in $B$ is given
by $C_o$ as before, we may assume that $B$ is realized as $B_{\btau}$ for some $\btau \in \mcM_1$. Consider a path
$\gamma:[0, 1] \to \mcM_1$ so that 
(i) $\gamma(0) = \btau, \gamma(1) = \btau_o$,
(ii) $B_{\gamma(s)} \in \mcR(\Cmb_{111})$ for $\forall s \in [0,1]$ and 
(iii) $t_{i_1} \rightsquigarrow -2$, $t_{i_2} \rightsquigarrow -1$, $t_{j_1} \rightsquigarrow 1$ and
$t_{j_2} \rightsquigarrow 2$. Then $B = B_{\btau}$ is deformed to $B_{\btau_o}$ or $B'_{\btau_o}$.
As we remark as above, $B_{\btau_o}$ is transformed to $B'_{\btau_o}$.
This shows that $B$ is continuously deformed to $B_{\btau_o}$ in $\mcR(\Cmb_{111})$.
Hence $\mcR(\Cmb_{111})$ is connected. 

\subsection{\texorpdfstring{$\Cmb_{121}$}{}}\label{subsec:121}
For $\btau = (\lambda,\bt) \in \mcM_2$, we define two elements $B_{\btau}, B'_{\btau}$ in $\mcR(\Cmb_{121})$
\[
 B_{\btau} =  C_o + L_{12, \btau} + L_{34, \btau} + D_{\btau} + L_{0, \btau}, \, 
 B'_{\btau} = C_o + L_{12, \btau} + L_{34, \btau} + D_{\btau} + L'_{0, \btau}.
\]
Here $D_{\btau}$ is the conic given by $c_{\btau}(t, x) = 0$.
Put $\ba = (-2, -1, 1, 2) \in \mcM$ and choose $\btau_o = ( \lambda_o,\ba) \in \mcM_2$ so that both of 
\[
 B_{\btau_o} =  C_o + L_{12, \btau_o} + L_{34, \btau_o} + D_{\btau_o} + L_{0, \btau_o}, \,
 B'_{\btau_o} = C_o + L_{12, \btau_o} + L_{34, \btau_o} + D_{\btau_o} + L'_{0, \btau_o}
\]
are in $\mcR(\Cmb_{121})$.
Now choose $B\in \mcR(\Cmb_{121})$ arbitrary.
By taking suitable coordinates of $\PP^2$ so that the conic in $B$ is given by $C_o$ as before, we may assume that $B$ is realized as $B_{\btau}$ for some $\btau \in \mcM_2$.
Consider a path
$\gamma:[0, 1] \to \mcM_2$ so that
(i) $\gamma(0) = \btau, 
\gamma(1) = \btau_o$,
(ii) $B_{\gamma(s)} \in \mcR(\Cmb_{121})$ for $\forall s \in [0,1]$ and 
(iii) $t_1 \rightsquigarrow -2$, $t_2 \rightsquigarrow -1$, $t_3 \rightsquigarrow 1$ and $t_4 \rightsquigarrow 2$. 
Then $B = B_{\btau}$ is deformed to $B_{\btau_o}$ or $B'_{\btau_o}$.
As $D_{\btau_o}$ is invariant under $(t, x) \mapsto (-t, x)$, $B_{\btau_o}$ is transformed to $B'_{\btau_o}$.
This shows that $B$ is continuously deformed to $B_{\btau_o}$ in $\mcR(\Cmb_{121})$.
Hence $\mcR(\Cmb_{121})$ is connected.

\subsection{\texorpdfstring{$\Cmb_{122}$}{}}\label{subsec:122}
Let $B$ be an arbitrary element in $\mcR(\Cmb_{122})$. By choosing coordinates of $\PP^2$ so that the
conic in $\mcP_1$ is given by $C_o$ and $D$ is tangent to $C_o$ at $z_o$, we may assume that $B$ is deformed to an element in $\mcR(\Cmb_{122})$ of the form
\[
B_{\btau} = C_o + L_{12, \btau} + L_{34, \btau} + L_{i_1j_1, \btau} + D_{i_2j_2, \btau},
\]
for some $\btau \in \mcM_1$, where $\{i_1, i_2\} = \{1,2\}, \{j_1, j_2\} = \{3, 4\}$.
Take $\ba = (-2, -1, 1,2)$, $\btau_o = (0,\ba) \in \mcM_1$ and consider
\[
B_{\btau_o} = C_o + L_{12, \btau_o} + L_{34, \btau_o} + L_{14, \btau_o} + D_{23, \btau_o},
\]
where $D_{23, \btau_o}$ is given by $x - 3t^2 +2 = 0$.
Now consider a path $\gamma: [0,1] \to \mcM_1$  so that (i) $\gamma(0) = \btau, \gamma(1) = \btau_o$, (ii) $B_{\gamma(s)} \in \mcR(\Cmb_{122})$ for $\forall s \in [0,1]$ and (iii) $t_1 \rightsquigarrow -2$, $t_2 \rightsquigarrow -1$, $t_3 \rightsquigarrow 1$ and $t_4 \rightsquigarrow 2$. 
Then we infer that $B$ is deformed to $B_{\btau_o}$. Hence $\mcR(\Cmb_{122})$ is connected.

\subsection{\texorpdfstring{$\Cmb_{123}$}{}}\label{subsec:123}

As for notation and terminology of this
subsection about elliptic surfaces, we use
those in Section~\ref{sec:preliminaries}.

We first show that there exists a Zariski pair $(B_1, B_2)$ for the combinatorics
$\Cmb_{123}$. 
Let $\mcC\mcL_{123} := \mcP_1\sqcup \mcA_2$ (
$\mcP_1 = \{C, L_1, L_2\}, \mcA_2=\{D, M\}$) be
a CL-arrangement with $\Cmb_{123}$. Let $\mcQ:= B_{\mcP_1}$ and
choose the tangent point between $C$ and $D$ as $z_o$. We assume
that $D$ is tangent to $L_1$ and $L_1 = L_{12}, L_2 = L_{34}$.
Let $S_{\mcQ, z_o}$ be the rational elliptic surface as before.  $D$ and $M$ give rise to a conic  point $P_D$ and
a line point $P_M$. By Table 1, we have 

\begin{eqnarray*}
    P_D & = & [\pm 1](P_{13}\dot{+}P_{14})\quad \mbox{or}\quad [\pm 1](P_{13}\dot{-}P_{14})\\
    P_M & = & [\pm 1](P_{13}\dot{+}P_{14}\dot{+}P_{34})
    \quad \mbox{or} \quad [\pm 1](P_{13}\dot{-}P_{14}\dot{+}P_{34})
\end{eqnarray*}
 Our tool to distinguish the embedded topology of CL-arrangement with $\Cmb_{123}$ is so called
 {\it the splitting types} introduced in \cite{bannai16} as follows:

\begin{defin}[{\cite{bannai16}}]
Let $\phi:X\to\PP^2$ be a double cover branched at a plane curve $\mcB$, and
let $D_1, D_2\subset\PP^2$ be two irreducible curves such that $\phi^\ast D_i$ are reducible and $\phi^\ast D_i=D_i^++D_i^-$.
For integers $m_1\leq m_2$, we say that the triple $(D_1, D_2;\mcB)$ has a {\it splitting type} $(m_1, m_2)$ if for a suitable choice of labels $D_1^+\cdot D_2^+=m_1$ and $D_1^+\cdot D_2^-=m_2$.
\end{defin}

The following proposition enables us to distinguish the embedded topology of plane curves by the splitting type.

\begin{prop}[{\cite[Proposition~2.5]{bannai16}}]
Let $\phi_i:X_i\to\PP^2$ $(i=1,2)$ be two double covers branched along plane curves $\mcB_i$, respectively.
For each $i=1,2$, let $D_{i1}$ and $D_{i2}$ be two irreducible plane curves such that $\phi_i^\ast D_{ij}$ are reducible and $\phi_i^\ast D_{ij}=D_{ij}^++D_{ij}^-$.
Suppose that $D_{i1}\cap D_{i2}\cap \mcB_i=\emptyset$, $D_{i1}$ and $D_{i2}$ intersect transversally, and that $(D_{11},D_{12};\mcB_1)$ and $(D_{21}, D_{22};\mcB_2)$ have distinct splitting types.
Then there is no homeomorphism $h:\PP^2\to\PP^2$ such that $h(\mcB_1)=\mcB_2$ and $\{h(D_{11}), h(D_{12})\}=\{D_{21}, D_{22}\}$.
\end{prop}

Under these setting, 
we have the following lemma:

\begin{lem}\label{lem:cmb123-split}{
$(D, M; \mcQ) = (0,2)$ if and only if $P_D\dot{+}P_M\dot{+}P_{34} = O$ with a suitable choice
of $P_D$ and $P_M$.
}
\end{lem}

\proof Let $s_D$ and $s_M$ be the sections corresponding to
$P_D$ and $P_M$, respectively. By \cite[Lemma~2.3]{BKMT22},
\[
s_D\cdot s_M = - \langle P_D, P_M\rangle + 1.
\]
Hence  $(D, M;\mcQ) = (0, 2)$ if and only if  $(\langle P_D, P_M\rangle, \langle P_D, [-1]P_M\rangle) = (1, -1)$ or $(-1,1)$. Now our statement follows from the following table:

\medskip
\begin{center}
    \begin{tabular}{|c|c|c|c|} \hline
      $P_D$ & $P_M$  & $\langle P_D, P_M\rangle$ & $s_D\cdot s_M$ \\ \hline
      $P_{13}\dot{\pm}P_{14}$ & $P_{13}\dot{
      \pm}P_{14}\dot{+}P_{34}$& $1$ & $0$ \\
       $P_{13}\dot{\pm}P_{14}$ & $[-1](P_{13}\dot{\pm}P_{14}\dot{+}P_{34})$ & $-1$ & $2$ \\
        $P_{13}\dot{\pm}P_{14}$ & $P_{13}\dot{\mp}P_{14}\dot{+}P_{34}$ & $0$ & $1$ \\
         $P_{13}\dot{\pm}P_{14}$ & $[-1](P_{13}\dot{\mp}P_{14}\dot{+}P_{34})$ & $0$ & $1$ \\ \hline
    \end{tabular}
    
\medskip
    Table 2
\end{center}


\qed

For $\mcC\mcL_{123}$, we can also take $\{D, M, L_{34}\}$ (resp.$\{C, L_{12}\}$) as $\mcP_1$ (resp.$\mcA_2$). 
Put $\mcQ' = D + M + L_{34}$.
Then we can also consider $(C, L_{12};\mcQ')$ and next lemma holds.

\begin{lem}\label{lem:cmb123-split2}{$(C, L_{12};\mcQ') = (0,2)$ if and only if 
$(D, M:\mcQ) = (0,2)$.
}
\end{lem}

\proof We choose homogeneous coordinates of $\PP^2$ as before. If $(D, M;\mcQ) = (0,2)$,
then we may assume that $P_D\dot{+}P_M\dot{+}P_{34} = O$.
Put $P_D = (x_{P_D}, y_{P_D}), P_M = (x_{P_M}, y_{P_M})$. Since the $x$-coordinate
of $P_D$ and $P_M$ give defining equations of $D$ and $M$, respectively,  we may assume
that $x_{P_D}, x_{P_M} \in \CC[t]$, $\deg x_{P_D} = 2, \deg x_{P_M} = 1$ and there exist
$mx + n \in \CC(t)[x]$ such that three points $P_D, P_M$ and $P_{34}$ are on the
line $y = mx + n$ in ${\mathbb A}^2_{\CC(t)}$. Put 
\[f_{\mcQ', z_o} = 
(x - x_{P_M})(x -x_{P_D})(x - (t_3+t_4)t - t_3t_4).
\]
Then we have
\[
f_{\btau}(t, x) - (mx +n)^2 = f_{\mcQ', z_o}, \qquad \btau=(0, \bt) \in \mcM_1. 
\]

Now consider a rational elliptic surface $S_{\mcQ', z_o}$ whose Weierstrass equation of
$E_{\mcQ', z_o}$ is given by $y^2 = f_{\mcQ', z_o}$.
From the above relation, three points $R_1, R_2$ and $R_3$ given by
\[
R_1:= (t^2, \sqrt{-1}(mt^2 + n)), \,\,   R_2:= (x_{P_{13}}, \sqrt{-1}(mx_{P_{13}}+ n)), \,\,  R_3:=(x_{P_{34}}, \sqrt{-1}(mx_{P_{34}}+ n)),
\]
where $x_{P_{ij}} = (t_i + t_j)t - t_it_j$, are on $y = \sqrt{-1}(mx + n)$. Hence 
$R_1\dot{+}R_2\dot{+}R_3 = O$. By Lemma~\ref{lem:cmb123-split}, $(C, L_{12};\mcQ') = (0,2)$. The converse statement follows by the same argument.
\qed

Now put
 \[
 B_1 = \mcQ + D + L_0, \,\, B_2 = \mcQ + D + L'_0
 \]
 where $D = \tilde{f}_{\mcQ, z_0}(s_{P_{13}\dot{+}P_{14}})$, 
 $L_0 = \tilde{f}_{\mcQ, z_o}(s_{P_{13}\dot{+}P_{14}\dot{+}P_{34}})$ and
 $L'_0 = \tilde{f}_{\mcQ, z_o}(s_{P_{13}\dot{-}P_{14}\dot{+}P_{34}})$. 
  Then we have

\begin{prop}\label{prop:zpair123}{$(B_1, B_2)$ is a Zariski pair.
}
\end{prop}

\proof Suppose that there exists a homeomorphism $h: (\PP^2, B_1) \to (\PP^2, B_2)$.
Then either $h(\mcQ) = \mcQ$ or $h(\mcQ') = \mcQ$ holds. Since
$(D, L_0;\mcQ) = (C, L_1;\mcQ') = (0, 2)$, $(D, L'_0;\mcQ) = (1, 1)$,
both cases are impossible by \cite[Proposition 2.5]{bannai16}.
\qed

\begin{rem} The $mx + n$ in our proof of Lemma~\ref{lem:cmb123-split} is in
$\CC[t,x]$ and its degree is $2$ as $f_{\btau}(t, x) - (mx +n)^2 = f_{\mcQ', z_o}$. Since $p_3, p_4$
and $p_0$ are on both $L_{34}$ and the conic $\tilde{C}$ given by $mx + n =0$ in the $(t, x)$-plane, we see that $\tilde C$ contains $L_{34}$. Hence we infer that
the three tangent point between $D+M$ and $C + L_{12}$ are collinear.
    
\end{rem}

We here give an explicit example of a Zariski pair for
$\Cmb_{123}$. We here keep previous notation.

\begin{ex}\label{ex:123}{\rm Let $\mcQ_{\btau_o}$ be a plane quartic given by $f_{\btau_o} = 0$ as before
where $\ba = (-2, -1, 1,2)$,$\btau_o = (0,\ba) \in \mcM_1$.
Let $S_{\mcQ_{\btau_o}, z_o}$ be the rational elliptic surface given by the Weierstrass equation $y^2 = f_{\mcQ_{\btau_o}, z_o}$ and $z_o = [0,1,0]$.
In this case, we have 
\[
P_{13}= (-t+2, 2\sqrt{2}(t-1)(t+2)), \quad P_{14} = (4, 3(t-2)(t+2))
\]
and $P_{D_{\btau_o}}:= P_{13}\dot{+}P_{14}$,  $P_{D'_{\btau_o}}:= P_{13}\dot{-}P_{14}$, $P_{M_{\btau_o}}:= P_D\dot{+}P_{34}$ and
$P_{M'_{\btau_o}}:= P_{D'_{\btau_o}}\dot{+}P_{34}$
as follows:
\begin{eqnarray*}
  P_{D_{\btau_o}} 
& = & (x_{D_{\btau_o}}, y_{D_{\btau_o}}) \\
   x_{D_{\btau_o}} &  = & (-12t^2 + 36t - 24)\sqrt{2} + 18t^2 - 51t + 34, \\
   y_{D_{\btau_o}} & = & 6(t - 2)(12\sqrt{2}t - 17\sqrt{2} - 17t + 24)(t - 1), \\
 P_{D'_{\btau_o}} & = & (x_{D'_{\btau_o}}, y_{D'_{\btau_o}}),  \\
 x_{D'_{\btau_o}} & = & (12t^2 - 36t + 24)\sqrt{2} + 18t^2 - 51t + 34 \\
 y_{D'_{\btau_o}} & = & 6(t - 2)(12\sqrt{2}t - 17\sqrt{2} + 17t - 24)(t - 1) 
 \end{eqnarray*}
 \begin{eqnarray*}
 P_{M_{\btau_o}} & = & (x_{M_{\btau_o}}, y_{M_{\btau_o}}) = (-2\sqrt{2}t -2,  -t^2 - \sqrt{2}t)\\
 P_{M'_{\btau_o}} & = & (x_{M'_{\btau_o}}, y_{M'_{\btau_o}})= (2\sqrt{2}t - 2,  t^2 - \sqrt{2}t) 
\end{eqnarray*}
Note that lines given by $x - x_{M_{\btau_o}}$ and $x - x_{M'_{\btau_o}}$ are
coincide with $L_{0, \btau_o}$ and $L'_{0, \btau_o}$, respectively.
Now put
\[
B_{1, \btau_o} = \mcQ_{\btau_o} + D_{\btau_o} + L_{0, \btau_o}, \quad
B_{2, \btau_o} = \mcQ_{\btau_o} + D_{\btau_o} + L'_{0, \btau_o}, 
\]
where $D_{\btau_o}$ is a conic of type $D(1,2)$ given by $ x - x_D = 0$.
We can easily check that
$B_{1, \btau_o}, B_{2, \btau_o} \in \mcR(\Cmb_{123})$  and $(B_{1, \btau_o}, B_{2, \btau_o})$ is a Zariski pair by
Proposition~\ref{prop:zpair123}.
}
\end{ex}

We now go on to study connected components of $\mcR(\Cmb_{123})$. 

\begin{prop}{Any element $B \in \mcR(\Cmb_{123})$ is deformed to either $B_{1, \btau_o}$ or
$B_{2, \btau_o}$ in Example~\ref{ex:123}, i.e., $\mcR(\Cmb_{123})$ has just
two connected components. 
}
\end{prop}

\proof By Example~\ref{ex:123}, $\mcR(\Cmb_{123})$ has at least two connected components.
Let $B$ be an element in $\mcR(\Cmb_{123})$. We show that $B$ is continuously deformed to $B_{1, \btau_o}$ or $B_{2, \btau_o}$ in Example~\ref{ex:123}. By taking homogeneous coordinates
suitably, we may assume that $B$ is given of the form
\[
B =  B_{\btau} = \mcQ_{\btau} + D_{\btau} + M_{\btau}, \quad \mcQ_{\btau} = C_o + L_{12,\btau} + L_{34, \btau}
\]
for some $\btau=(0,\bt) \in \mcM_1$ and $D_{\btau}$ and $M_{\btau}$ are the conic and line described in Section~\ref{subsec:explicit-construction}.
We may also assume that $D_{\btau}$ passes through $p_3$ and $p_4$ and tangent to $L_{12, \btau}$. Now consider a path $\gamma:[0,1] \to \mcM_1$ such that
(i) $\gamma(0) = \btau$, $\gamma(1) = \btau_o$,
(ii) $B_{\gamma(s)} \in \mcR(\Cmb_{123})$ for $\forall s \in [0,1]$ and
(iii)$t_1 \rightsquigarrow -2, t_2 \rightsquigarrow -1,
t_3 \rightsquigarrow 1$ and
$t_4 \rightsquigarrow 2$.
Then $B_{\gamma(0)} = B$ and 
\[
B_{\gamma(1)} = \mcQ_{\btau_o} + D_{12, \gamma(1)} + M_{\gamma(1)},
\]
where
\[
    D_{\gamma(1)}  =  D_{\btau_o}\,\, \mbox{or}\,\, D'_{\btau_o}, \quad 
    M_{\gamma(1)}  =  L_{0, \btau_o}\,\,  \mbox{or}\,\,   L'_{0,\btau_o},
\]
where $D'_{\btau_o}$ is the conic given by $x - x_{D'}=0$.

Case (i): $D_{\gamma(1)}= D_{\btau_o}$.
In this case, $B_{\gamma(1)}$ is either $B_{1, \btau_o}$ or $B_{2, \btau_o}$.

Case (ii): $D_{\gamma(1)} = D'_{\btau_o}$. 
In this case,
$B_{\gamma(1)}$ is either $\mcQ_{\btau_o} + D'_{\btau} + L_{0, \btau_o}$ or  $\mcQ_{\btau_o} + D'_{\btau} + L'_{0, \btau_o}$.
Consider families of lines and parabolas
as follows:
\begin{eqnarray*}
    L_{u_1u_2}&&: x - (u_1 + u_2)t + u_1u_2 = 0, \, (u_1, u_2) \in \CC^2, u_1\neq u_2, \\
    D_{\mu}&&: x - \mu t^2 - (3- 3\mu)t - 2\mu + 2 = 0, \, \mu \in \CC^{\times}.
\end{eqnarray*}
Namely, $L_{u_1u_2}$ is a line intersecting $C_o$ at $(u_1, u_1^2)$ and $(u_2, u^2_2)$ and
$D_{\mu}$ is a parabola passing $(1,1)$ and $(2, 4)$. It is easily checked that
the condition for $L_{u_1u_2}$ and $D_{\mu}$ to be tangent is
that $(u_1, u_2, \mu)$ satisfies 
\[
(\ast) \quad \mu^2 - 4u_1u_2\mu + 6(\mu-1)(u_1 + u_2) + (u_1 + u_2)^2 - 10\mu + 9 = 0.
\]
Note that  the surface given by $(\ast)$ in the
$(u_1, u_2, \mu)$-space is irreducible and connected.
Now consider a path $\bar{\gamma}: [0,1] \to \mcM_1\times \CC^{\times}, 
\bar{\gamma}(s) = (0, u_1(s), u_2(s), 1, 2, \mu(s))$ such that (i) $(u_1(s), u_2(s), \mu(s))$
satisfies $(\ast)$ and (ii) $\bar{\gamma}(0) = (0, -2, -1, 1,2, 18+12\sqrt{2})$ and $\bar{\gamma}(1) = (0, -2, -1, 1,2, 18- 12\sqrt{2})$.
Since
(i) $D_{\mu(s)}$ is tangent to
$L_{u_1u_2}$ and (ii) the line $M_{\bar{\gamma}(s)}$ is determined by $L_{u_1u_2}\cap L_{34, \ba}$ and the initial line $M_{\bar{\gamma}(0)}$, 
we infer that there exists a continuous family $B_{\bar{\gamma}(s)}$ ($0 \le s \le 1$) in $\mcR(\Cmb_{123})$ such that $B_{\bar{\gamma}(0)} = B_{\gamma(1)}$ and $B_{\bar{\gamma}(1)} = B_{1, \btau_o}$ or $B_{2, \btau_o}$.
Thus our statement follows.
\endproof

\subsection{\texorpdfstring{$\Cmb_{124}$}{}}\label{subsec:124}
In \cite{tokunaga14}, we have seen that there  exists a Zariski pair for $\Cmb_{124}$. Hence $\mcR(\Cmb_{124})$ has at least two connected components.
In this subsection, we will show that there exist only two components.
Let us start with the following example.

\begin{ex}\label{ex:cmb124} 
Let $\btau_o=(0,\ba) \in \mcM_1$ and $\mcQ_{\btau_o}$ be as before.
We label $p_1, p_2, p_3$ and $p_4$ in the same way. Namely lines contained in $\mcQ_{\btau_o}$ are $L_{12,\btau_o}$ and $L_{34, \btau_o}$. In this case, we have
\[
[2]P_{13} = \left (\frac 98 t^2, \frac{\sqrt{2}}{32}(-9t^3 + 16t) \right ), \,\,
[2]P_{14} = \left (t^2 + \frac 14, \frac 12 t^2 - \frac 98 \right ).
\]
Now put
\[
\begin{array}{lc}
D_{\btau_o}: \tilde{f}_{\mcQ_{\btau_o}, z_o}(s_{[2]P_{13}}) = x - \frac 98 t^2 = 0, &
L_{13, \btau_o}: \tilde{f}_{\mcQ_{\btau_o}, z_o}(s_{P_{13}}) = x + t - 2 = 0, \\
L_{14, \btau_o}: \tilde{f}_{\mcQ_{\btau_o}, z_o}(s_{P_{14}}) = x -4 = 0. & 
\end{array}
\]
Now define $B_1$ and $B_2$ to be
\[
B_{1, \btau_o} = \mcQ_{\btau_o} + D_{\btau_o} + L_{13, \btau_o}, \,\, 
B_{2, \btau_o} =  \mcQ_{\btau_o} + D_{\btau_o} + L_{14, \btau_o}.
\]
Then by \cite[Theorem 5]{tokunaga14}, $(B_{1,\btau_o}, B_{2,\btau_o})$ is a Zariski pair.
\end{ex}

Now we show 

\begin{prop}\label{prop:cmb124}{Let $B$ be an arbitrary member in $\mcR(\Cmb_{124})$.
Then $B$ is continuously deformed to either $B_{1,\btau_o}$ or $B_{2,\btau_o}$ in Example~\ref{ex:cmb124}.
In particular, $\mcR(\Cmb_{124})$ has just two connected components.
} 
\end{prop}
\proof After taking a suitable coordinate change, we may assume that
$B$ is given of the form 
\[
B =  B_{\btau} = \mcQ_{\btau} + D_{\btau} +  L_{13, \btau},
\]
for some $\btau = (0, \bt) \in \mcM_1$.
Here $D_{\btau}$ is either $\tilde{f}_{\mcQ_{\btau}, z_o}(s_{[2]P_{13, \btau}})$ or $\tilde{f}_{\mcQ_{\btau}, z_o}(s_{[2]P_{14, \btau}})$

\underline{Case $D_{\btau}=\tilde{f}_{\mcQ_{\btau}, z_o}(s_{[2]P_{13, \btau}})$.}
Consider a path $\gamma: [0,1] \to \mcM_1$ such that
(i) $\gamma(0) = \btau$, $\gamma(1) = \btau_o$,
(ii) $B_{\gamma(s)} \in \mcR(\Cmb_{124})$ for $\forall s \in [0,1]$ and
(iii) $t_1\rightsquigarrow -2$, $t_2\rightsquigarrow -1$, $t_3\rightsquigarrow 1$ and $t_4\rightsquigarrow 2$.
This shows that $B$ is continuously deformed to $B_{1,\btau_o}$ with keeping the combinatorics.

\underline{Case $D_{\btau}=\tilde{f}_{\mcQ_{\btau}, z_o}(s_{[2]P_{14, \btau}})$.}
Consider a path $\gamma: [0,1] \to \mcM_1$ such that
(i) $\gamma(0) = \btau$, $\gamma(1) = \btau_o$,
(ii) $B_{\gamma(s)} \in \mcR(\Cmb_{124})$ for $\forall s \in [0,1]$ and
(iii) $t_1\rightsquigarrow -2$, $t_2\rightsquigarrow -1$, $t_3\rightsquigarrow 2$ and $t_4\rightsquigarrow 1$.
Then $L_{13, \btau}$ (resp. $L_{14, \btau}) $ is deformed to $L_{14, \btau_o}$
(resp. $L_{13, \btau_o}$) and
$D_{\btau}$ is deformed to $\tilde{f}_{\mcQ_{\btau}, z_o}(s_{[2]P_{13, \btau_o}})$ accordingly.
Hence
$B$ is continuously deformed to $B_{2,\btau_o}$ with keeping the combinatorics.
 \endproof

\subsection{\texorpdfstring{$\Cmb_{125}$}{}}\label{subsec:125}
Choose $B \in \mcR(\Cmb_{125})$ arbitrary. By taking  appropriate coordinates of 
$\PP^2$, we may assume that $C_1 = C_o$, $D$ is tangent to $C_o$ at $z_o =[0,1,0]$ and
there exists $\btau = (0,\bt) \in \mcM_1$ such that $B$ is of the form
\[
B_{\btau}= C_o + L_{12,\btau} + L_{34, \btau} + D_{\btau} +  L_{0, \btau}.
\]
By Table 1, $D_{\btau}$ is given by the image of $s_{[2]P_{ij}}$ under ${\tilde f}_{\mcQ_{\btau}, z_o}$. Take $\ba=(-2, -1, 1, 2)$, $\btau_o = (\lambda_1,\ba) \in \mcM_1$ and consider an element of $\mcR(\Cmb_{125})$ given by
\[
B_{\btau_o} = C_o + L_{12, \btau_o} + L_{34, \btau_o} + D_{\btau_o} + L_{0, \btau_o},
\]
where $D_{\btau_o}$ is given by $[2]P_{23}$. By \cite[Example 5.2]{tokunaga14}, 
$D_{\btau_o}$ is given by $x - t^2 - 1/4 = 0$.
Now we choose a path $\gamma: [0,1] \to \mcM_1$ such that (i) $\gamma(0) = \btau, \gamma(1) = \btau_o$, (ii) $B_{\gamma(s)} \in \mcR(\Cmb_{125})$ for $\forall s \in [0,1]$ and (iii) $t_i \rightsquigarrow -1, t_j \rightsquigarrow 1$. By the deformation along $\gamma$, $L_{ij}$ is deformed to
$L_{23}$. Hence $D_{\btau}$ is deformed to $D_{\btau_o}$.
If $L_{0, \btau}$ is deformed to $L_{0, \btau_o}$, we see that $B$ is deformed to $B_{\btau_o}$.
If $L_{0, \btau}$ is deformed to $L'_{0, \btau_o}$, we then apply the transformation $(t, x) \mapsto (-t, x)$ and we see that $B$ is deformed to $B_{\btau_o}$.
Thus $\mcR(\Cmb_{125})$ is connected.

\subsection{\texorpdfstring{$\Cmb_{211}$}{}}\label{subsec:211}
Let us start with the following remark.

\begin{rem}\label{rem:bitangents}{Let  $B_{\mcP_2}$ be a quartic given by a conic arrangement $\mcP_2$. It is known that there exist
four bitangent lines for $B_{\mcP_2}$. When we deform
conics in $\mcP_2$ continuously, these bitangents are
also deformed along with conics. Note that this observation
follows by considering dual curves of the conics in $\mcP_2$. We make use of 
this observation repeatedly in the rest of this article.
}
\end{rem}

Consider two conics $C_{o1}$ and $C_{o2}$ given by 
\[
C_{o1}:  t^2+x^2+tx-\frac{27}{4}=0, \quad C_{o2}: t^2+x^2-tx-\frac{27}{4}=0.
\]
We write $C_{o1}\cap C_{o2}$ by   $\bp=(p_1,p_2,p_3,p_4)$ whose affine
coordinate is given by
\[
p_1 = \left (0,\frac{3}{2}\sqrt{3}\right ), \quad
p_2 = \left (\frac{3}{2}\sqrt{3},0 \right ), \quad
p_3 = \left (0,-\frac{3}{2}\sqrt{3}\right ), \quad
p_4 = \left (-\frac{3}{2}\sqrt{3},0 \right ). \quad
\]
The bitangent lines of $C_{o1} + C_{o2}$ are
\[
L_{b1,\bp}: t = 3, \,\,
L_{b2,\bp}: t = -3, \,\,
L_{b3,\bp}: x = 3, \, \,
L_{b4,\bp}: x = -3.
\]
Now put
\[
B_{oi} := C_{o1} + C_{o2} + L_{13, \bp} + L_{24, \bp} + L_{bi,\bp}, \,\, i = 1, 2, 3, 4.
\]
Then $B_{oi} \in \mcR(\Cmb_{211})$ and all of them are transformed by some
projective transformation each other.

Hence it is enough to show an arbitrary element $B \in \mcR(\Cmb_{211})$ can be continuously deformed to $B_{oi} \in \mcR(\Cmb_{211})$ for some $i$.

We may assume that $B$ is given in the follow form:
\[
B_{\btau} = \mcQ_{\btau} + L_{13, \btau} + L_{24, \btau} + L_{b_1, \btau},
\]
where $\mcQ_{\btau} = C_o + C_{\btau}$ for some $\btau=(\lambda,\bt) \in \mcM_2$. Let $\phi: \PP^2 \to \PP^2$ be a projective transformation such that $\phi(C_o) = C_{o1}$.
Then there exists $\btau_{\bc}=(\lambda_{\bc},\bc) \in \mcM_2$ such that $\phi(C_{\btau_{\bc}}) =C_{o2}$ and points in $C_o\cap C_{\btau_{\bc}}$ are labeled so that $L_{ij, \btau_{\bc}} = L_{ij, \bp}$ holds.
 Now we choose a path $\gamma: [0,1] \to \mcM_2$
such that
(i) $\gamma(0) = \btau, \gamma(1) = \btau_{\bc}$,
(ii) $B_{\gamma(s)} \in \mcR(\Cmb_{211})$ for $\forall s \in [0,1]$ and
(iii) $t_i \rightsquigarrow c_i$ $(i =1, 2, 3, 4)$.
We see that $B$ can be continuously deformed along $\gamma$ in $\mcR(\Cmb_{211})$ to $B_1:=C_o + C_{\btau_{\bc}} + L_{13, \btau_{\bc}}+ L_{24, \btau_{\bc}} + L_{b,\btau_{\bc}}$.
Here $L_{b,\btau_{\bc}}$ denotes a bitangnet to $C_o + C_{\btau_{\bc}}$.
As $\phi(B_1) = B_{oi}$ for some $i$, we infer that
$B$ is continuously deformed to $B_{oi}$ and that
$\mcR(\Cmb_{211})$ is connected.

\subsection{\texorpdfstring{$\Cmb_{212}$}{}}
We first show that there exists a Zariski pair for $\Cmb_{212}$.
Let $\mcQ_{\btau} = C_o + C_{\btau}$ and $B = \mcQ_{\btau} +L_{ij} + L_{bk} + L_{bl} \in \mcR(\Cmb_{212})$. 
Choose $z_o \in C_o$ so that the tangent line at $z_o$ meets $C_{\btau}$ at two distinct points.
Let $\varphi_{\mcQ_{\btau}, z_o} : S_{\mcQ_{\btau}, z_o} \to \PP^1$ and $\tilde{f}_{\mcQ_{\btau}, z_o}: S_{\mcQ_{\btau}, z_o} \to \PP^2$ as in Subsection~\ref{subsec:elliptic_surface}. As we have seen in Table 1 or \cite[Section 3.2]{BKMT22}, if we put
\[
\begin{array}{cc}
    Q_1:= P_{12} \dot{+} P_{13}\dot{+}P_{23}, &
    Q_2:= [-1]P_{12} \dot{+} P_{13}\dot{+}P_{23}, \\
    Q_3:= P_{12} \dot{+} [-1]P_{13}\dot{+}P_{23}, &
    Q_4:= P_{12} \dot{+} P_{13}\dot{+}[-1]P_{23}, \\
\end{array}
\]
then we may assume that $L_{bi}:= \tilde{f}_{\mcQ_{\btau}, z_o}(s_{Q_i})$ $(i = 1,2, 3, 4)$ are four bitangent lines of $\mcQ_{\btau}$.
Then by \cite[Theorems 3.2 and 3.3]{tokunaga14} and the argument in p.629-630 in \cite{tokunaga14}, we have the following proposition:
\begin{prop}\label{prop:212}
Let $p$ be an odd prime.
There exits a $D_{2p}$-cover of $\PP^2$ branched at $2\mcQ_{\btau} + p(L_{ij} + L_{bk} + L_{bl})$ if and only if the images of $P_{ij}, Q_k, Q_l$ in $E_{\mcQ_{\lambda}, z_o}$ are linearly dependent over $\ZZ/p\ZZ$.    
\end{prop}

By Proposition~\ref{prop:212}, we have
\begin{cor}\label{cor:212}{Let $B_{kl}:= \mcQ_{\btau} + L_{13} + L_{bk} + L_{bl}$.
Then $(B_{13},B_{kl})$ (resp.$(B_{24},B_{kl})$) is Zariski pairs where $(k,l) \neq (2,4)$ (resp.$(k,l) \neq (1,3)$).
}    
\end{cor}
\proof If a homeomorphism $h : (\PP^2, B_{13}) \to (\PP^2, B_{kl})$ exists,
it satisfies $h(\mcQ_{\btau}) = \mcQ_{\btau}$. Hence our statement follows from
Proposition~\ref{prop:212}.
\endproof

\begin{rem}
 {\rm We may use the connected number for $L_{13} + L_{bk} + L_{bl}$ to prove our statement. In fact,
 for example, the connected number is 2 for $(k, l) = (1, 3)$, while it is $1$ for $(k,l) = (1, 2)$. This shows $(B_{12}, B_{13})$ is a Zariski pair.
 As for connected numbers, see \cite{shirane18} for detail.}  
 \end{rem}
 
Let us now consider an explicit example.

\begin{ex}
\label{ex:cmb212a}
Let $\mcQ_{\btau_o} = C_o + C_{\btau_o}$ be a plane quartic given by $f_{\mcQ_{\btau_o}}=0$ where
$\btau_o = (\lambda_o,\ba) = (-10, -2, -1, 1,2)$.
Let $S_{\mcQ_{\btau_o}}$ be the rational elliptic surface given by the Weierstrass equation $y^2=f_{\mcQ_{\btau_o}}$ and $z_o=[0,1,0]$.
In this case, we have
\begin{eqnarray*}
P_{12} &=& \left( -3t-2, -i \, \sqrt{10} t^{2} - 3 i \, \sqrt{10} t - 2 i \, \sqrt{10} \right), \\
P_{13} &=& \left( -t+2, -i \, \sqrt{2} t^{2} - i \, \sqrt{2} t + 2 i \, \sqrt{2} \right), \\
P_{23} &=& \left( 1,  -i \, t^{2} +  i \right).
\end{eqnarray*}

Under these setting,  $P_{L_{b1}}:=P_{12} \dot{+}P_{13} \dot{+}P_{23}$, $P_{L_{b2}}:=P_{12} \dot{+}P_{13} \dot{-}P_{23}$, $P_{L_{b3}}:=P_{12} \dot{-}P_{13} \dot{+}P_{23}$ and $P_{L_{b4}}:=[-1]P_{12} \dot{+}P_{13} \dot{+}P_{23}$ 
are given are as follows:
\begin{eqnarray*}
P_{L_{b1}} &=& \left( \sqrt{2}(\sqrt{5}+3)t-3\sqrt{5}-7,
 (2\sqrt{5} + 3)it^{2} - \frac{\sqrt{2}}{2}(15\sqrt{5} + 29) it + 2(7 \sqrt{5} + 15) i  \right), \\
P_{L_{b2}} &=& \left( -\sqrt{2}(\sqrt{5}+3)t-3\sqrt{5}-7,
 (2\sqrt{5} + 3)it^{2} + \frac{\sqrt{2}}{2}(15\sqrt{5} + 29) it + 2(7 \sqrt{5} + 15) i\right), \\
P_{L_{b3}} &=& \left( \sqrt{2}(\sqrt{5}-3)t+3\sqrt{5}-7,
 -(2\sqrt{5} - 3)it^{2} - \frac{\sqrt{2}}{2}(15\sqrt{5} - 29) it - 2(7 \sqrt{5} - 15) i\right), \\
P_{L_{b4}} &=& \left( -\sqrt{2}(\sqrt{5}-3)t+3\sqrt{5}-7,
 -(2\sqrt{5} - 3)it^{2} + \frac{\sqrt{2}}{2}(15\sqrt{5} - 29) it - 2(7 \sqrt{5} - 15) i \right).
\end{eqnarray*}
Now put $L_{bi,\btau_o}:f_{\mcQ_{\btau_o},z_o}(s_{P_{bi}}) = 0$. Then we have
\[
\begin{array}{cc}
L_{b1,\btau_o} : x - \sqrt{2}(\sqrt{5}+3)t +3\sqrt{5} +7 =0, &
L_{b2,\btau_o} : x + \sqrt{2}(\sqrt{5}+3)t +3\sqrt{5} +7 =0, \\
L_{b3,\btau_o} : x - \sqrt{2}(\sqrt{5}-3)t -3\sqrt{5} +7 =0, &
L_{b4,\btau_o} : x + \sqrt{2}(\sqrt{5}-3)t -3\sqrt{5} +7 =0.
\end{array}
\]
We put
\[
B_{ij,\btau_o}=\mcQ_{\btau_o} + L_{13,\btau_o} + L_{bi,\btau_o} + L_{bj,\btau_o}, \quad i,j = 1,2,3,4, i\neq j.
\]
Then $(B_{13,\btau_o},B_{ij,\btau_o})$ (resp.$(B_{24,\btau_o},B_{ij,\btau_o})$) are Zariski pairs for $(i,j) \neq (2,4)$ (resp.$(i,j) \neq (1,3)$) by Corollary~\ref{cor:212}.
\end{ex}

We give another example of a CL-arrangement in $\Cmb_{212}$, which plays an important role to study the connectivity for $\mcR(\Cmb_{212})$.
\begin{ex}
\label{ex:cmb212p}
Let $C_{o1}$ and $C_{o2}$ be conics given by
\[
C_{o1}: t^2 + x^2 + tx -\frac{27}{4} = 0, \qquad
C_{o2}: 676t^2 + 764tx + 676x^2 - 4563 = 0.
\]
We write $C_{o1}\cap C_{o2}$ by $\bp=(p_1,p_2,p_3,p_4)$ whose affine coordinate is given by
\[
p_1=\left( \frac{3}{2}\sqrt{3},0 \right), \quad
p_2=\left( 0, -\frac{3}{2}\sqrt{3} \right), \quad
p_3=\left( -\frac{3}{2}\sqrt{3},0 \right), \quad
p_4=\left( 0, \frac{3}{2}\sqrt{3} \right).
\]
The bitangent lines of $C_{o1}+C_{o2}$ are
\[
L_{b1,\bp}:15t + 8x - 39=0,
L_{b2,\bp}:8t + 15x + 39=0,
L_{b3,\bp}:15t + 8x + 39=0,
L_{b4,\bp}:8t + 15x - 39=0.
\]
Now put
\[
B_{ij,\bp}=C_{o1} + C_{o2} + L_{13,\bp} + L_{bi,\bp} + L_{bj,\bp}, \quad i,j = 1,2,3,4.
\]
Then $B_{ij,\bp} \in \mcR(\Cmb_{212})$.
\end{ex}

Now we show
\begin{prop}
\label{prop:cmb212}
There exists a homeomorphism $h:\PP^2 \to \PP^2$ such that $h(B_{13,\bp})=B_{24,\bp}$.
\end{prop}
\proof
Let $\beta$ be a parameter and $C_{o2,\beta}$ be a conic defined by
\begin{eqnarray*}
C_{o2,\beta}&:&-27 \, \beta^{4} - 54 \, \beta^{3} + 4 \, {\left(\beta^{4} + 2 \, \beta^{3} + 3 \, \beta^{2} + 2 \, \beta + 1\right)} t^{2} + 4 \, {\left(2 \, \beta^{4} + 4 \, \beta^{3} - 6 \, \beta^{2} - 8 \, \beta - 1\right)} t x \\
&&+ 4 \, {\left(\beta^{4} + 2 \, \beta^{3} + 3 \, \beta^{2} + 2 \, \beta + 1\right)} x^{2} - 81 \, \beta^{2} - 54 \, \beta - 27 = 0.
\end{eqnarray*}
The conic $C_{o2,\beta}$ passes through $p_1,p_2,p_3,p_4$ and furthermore, $C_{o2,\beta}=C_{o2}$ for $\beta=-4,-\frac{5}{7},-\frac{2}{7},3$.
Note that $L_{13,\bp}$ is fixed since $p_i$ does not depend on the parameter $\beta$.
Also, $C_{o2,\beta}=C_{o1}$ if $(\beta^2+4\beta+1)(\beta^2-2\beta-2)=0$ and $C_{o2,\beta}$ has singular points if $(\beta^2+\beta+1)(2\beta^2+2\beta-1)(2\beta+1)=0$.

The three lines $L_{13,\bp},L_{bi,\bp,\beta},L_{bj,\bp,\beta}$ $(i,j \in \{1,2,3,4\})$ intersect at one points if $\beta=-2,-1,0,1$.
Now we have the following bitangent lines $L_{bi,\bp,\beta}$ of $C_{o1}+C_{o2,\beta}$ when $\beta(\beta-1)(\beta+1)(\beta+2)(\beta^2+4\beta+1)(\beta^2-2\beta-2)(\beta^2+\beta+1)(2\beta^2+2\beta-1)(2\beta+1)\neq 0$:
\begin{align*}
&L_{b1,\bp,\beta} : (\beta^2+2\beta)t + (\beta^2-1)x - (3\beta^2 + 3\beta + 3)=0 \\
&L_{b2,\bp,\beta} : (\beta^2-1)t + (\beta^2+2\beta)x + (3\beta^2 + 3\beta + 3)=0 \\
&L_{b3,\bp,\beta} : (\beta^2+2\beta)t + (\beta^2-1)x + (3\beta^2 + 3\beta + 3)=0 \\
&L_{b4,\bp,\beta} : (\beta^2-1)t + (\beta^2+2\beta)x - (3\beta^2 + 3\beta + 3)=0.
\end{align*}
For $\beta=-4,-\frac{5}{7},-\frac{2}{7},3$, we have the following table:
\begin{table}[H]
\centering
\begin{tabular}{|c||c|c|c|c|}
\hline
$\beta$ & $-4$ & $-\frac{5}{7}$ & $-\frac{2}{7}$ & $3$\\ \hline
$L_{b1,\bp,\beta}$ & $L_{b4,\bp}$ & $L_{b3,\bp}$ & $L_{b2,\bp}$ & $L_{b1,\bp}$ \\ \hline
$L_{b2,\bp,\beta}$ & $L_{b3,\bp}$ & $L_{b4,\bp}$ & $L_{b1,\bp}$ & $L_{b2,\bp}$ \\ \hline
$L_{b3,\bp,\beta}$ & $L_{b2,\bp}$ & $L_{b1,\bp}$ & $L_{b4,\bp}$ & $L_{b3,\bp}$ \\ \hline
$L_{b4,\bp,\beta}$ & $L_{b1,\bp}$ & $L_{b2,\bp}$ & $L_{b3,\bp}$ & $L_{b4,\bp}$ \\ \hline
\end{tabular}
\end{table}

By considering $C_{o1}+C_{o2,\beta}+L_{13,\bp}+L_{b1,\bp,\beta}+L_{b3,\bp,\beta}$ and $C_{o1}+C_{o2,\beta}+L_{13,\bp}+L_{b2,\bp,\beta}+L_{b4,\bp,\beta}$ for $\beta=-4,-\frac{5}{7},-\frac{2}{7},3$, we see that
\[
C_{o1}+C_{o2}+L_{13,\bp}+L_{b1,\bp}+L_{b3,\bp} \sim C_{o1}+C_{o2}+L_{13,\bp}+L_{b2,\bp}+L_{b4,\bp}
\]
since we can deform while avoiding the finite number of exceptional values of $\beta$ where the combinatorics becomes degenerated.
Hence our statement follows.
\qed

\begin{rem}
By the proof in the above Proposition, we see that there also exists a homeomorphism $h:\PP^2 \to \PP^2$ such that $h(B_{ij,\bp})=B_{kl,\bp}$ for $(i,j,k,l)=(1,2,3,4)$,  $(1,4,2,3)$ or $(1,3, 2, 4)$.
\end{rem}

\begin{cor}
\label{cor:cmb212}
There exists a homeomorphism $h:\PP^2 \to \PP^2$ such that $h(B_{12,\bp})=B_{14,\bp}$.
\end{cor}
\proof
We use the same example in Proposition \ref{prop:cmb212}.
We put 
\begin{eqnarray*}
B_{12,\bp,\beta}&:=&C_{o1}+C_{o2,\beta}+L_{13,\bp}+L_{b1,\bp,\beta}+L_{b2,\bp,\beta}, \\
B_{14,\bp,\beta}&:=&C_{o1}+C_{o2,\beta}+L_{13,\bp}+ L_{b1,\bp,\beta}+L_{b4,\bp,\beta}.
\end{eqnarray*}
By letting $\beta'=0$, we see that 
$C'_{o2}:=C_{o2,\beta'}$ is given by
\[
C'_{o2}: t^2 + x^2 - tx - \frac{27}{4} = 0.
\]
and the bitangent lines of $C_{o1}+C'_{o2}$ are
\[
L_{b1,\bp,\beta'}: x - 3 =0, \quad
L_{b2,\bp,\beta'}: t + 3 =0, \quad
L_{b3,\bp,\beta'}: x + 3 =0, \quad
L_{b4,\bp,\beta'}: t - 3 =0.
\]
Then $B_{12,\bp,\beta'},B_{14,\bp,\beta'} \in \mcR(\Cmb_{212})$ are transformed by $[T, X, Z] \mapsto [-T, X, Z]$.
Hence $B_{12,\bp}$ can be deformed to $B_{14,\bp}$. 
Hence our assertion follows.
\qed

We are now in position to prove the following proposition:
\begin{prop}\label{prop:212-final}{
Any element $B \in \mcR(\Cmb_{212})$ is deformed to either $B_{12,\bp}$ or $B_{13,\bp}$ in Example \ref{ex:cmb212p}, i.e., $\mcR(\Cmb_{212})$ have just two connected components.
}
\end{prop}
\proof
Our proof consists of two steps:
\begin{enumerate}
\item[(I)] Any element $B \in \mcR(\Cmb_{212})$ is deformed to $B_{ij,\btau_o}$ $(i,j \in \{1,2,3,4\}, i\neq j)$ in Example \ref{ex:cmb212a}.
\item[(II)] $B_{ij,\btau_o}$ is deformed to either $B_{12,\bp}$ or $B_{13,\bp}$ in Example \ref{ex:cmb212p}. 
\end{enumerate}
Since $B_{12, \bp}$ and $B_{13, \bp}$ belong to distinct connected 
components of $\mcR(\Cmb_{212})$, Steps (I) and (II) implies Proposition~\ref{prop:212-final}. 

\underline{Step (I):}
After taking a suitable coordinates change and labeling the intersection points
$C_1\cap C_2$, we may assume that $B$ is given as follows:

There exists $\btau \in \mcM_2$ such that
\[
B_{\btau} = \mcQ_{\btau} + L_{i_1i_2,\btau} + L_{b{j_1},\btau} + L_{b{j_2},\btau}, \quad i_1,i_2,j_1,j_2 \in \{1,2,3,4\}
\]
where $L_{b{j_1},\btau}$ and $L_{b{j_2},\btau}$ are given $\tilde{f}_{\mcQ_{\btau}, z_o}( Q_{j_1})$ and 
$\tilde{f}_{\mcQ_{\btau}, z_o}(Q_{j_2})$, respectively.

Now consider a path $\gamma:[0,1]\to \mcM_2$ such that 
(i) $\gamma(0)=\btau$, $\gamma(1)=\btau_o$
(ii) $B_{\gamma(s)} \in \mcR(\Cmb_{212})$ for $\forall s \in [0,1]$, and
(iii) $t_{1}\rightsquigarrow -2$, 
$t_{2}\rightsquigarrow -1$, $t_{3}\rightsquigarrow 1$,
$t_{4}\rightsquigarrow 2$.
Then $B$ is deformed to $B_{ij, \btau_o}$.
Hence we have the assertion in Step (I).

\medskip

\underline{Step (II):} Let $B_{ij, \btau_o}$ be the CL-arrangement as in Example~\ref{ex:cmb212a}.
By Corollary \ref{cor:212}, $\mcR(\Cmb_{212})$ has at least two connected components.
We here show that any $B_{ij, \btau_o}$ which has 6 possibilities can be continuously deformed to 
either $B_{12,\bp}$ or $B_{13,\bp}$.

Let $\phi : \PP^2 \to \PP^2$ be a projective transformation  such that $\phi(C_o) = C_{o1}$.
We choose $\bc=(c_1, c_2, c_3, c_4)$ and $\btau_{\bc}=(\lambda_{\bc}, \bc) \in \mcM_2$ such that $\phi(\mcQ_{\btau_{\bc}}) = C_{o1} + C_{o2}$.
Now we choose path $\gamma$  in $\mcM_2$ as in Step (I) such that $\gamma(0) = \btau_{\bc}$ and $\gamma(1) =\btau_o$.
Then we infer that $B_{ij, \btau_o}$ is continuously deformed to $B_{i_1j_1, \btau_{\bc}}$ in $\mcR(\Cmb_{212})$.
Since $\phi(B_{i_1j_1, \btau_{\bc}}) =  B_{i_2j_2,\bp}$ for some $i_2, j_2$, we see that 
$B_{ij, \btau_o}$ is continuously deformed to $B_{i_2j_2, \bp}$.
Now by Proposition \ref{prop:cmb212} and Corollary \ref{cor:cmb212}, $B_{ij,\bp}$ is deformed to either $B_{12,\bp}$ or $B_{13,\bp}$ and we have the assertion in Step (II).
\endproof

\subsection{\texorpdfstring{$\Cmb_{213}$}{}}\label{subsec:213}
We keep our notation in $\Cmb_{211}$. Let $B$ be an arbitrary element in $\mcR(\Cmb_{213})$ and
we may assume that $B$ is given in the form
\[
B_{\btau} = \mcQ_{\btau} + L_{b_1, \btau} + L_{b_2, \btau} + L_{b_3, \btau}
\]
for some $\btau =(\lambda,\bt) \in \mcM_2$. 
In other words,  $B$ is determined by the remaining bitangent $L_{b_4, \btau}$.
Hence we infer that it is enough to show that $\mcQ_{\btau} + L_{b_4, \btau}$ can be continuously deformed to $C_{o1} + C_{o2} + L_{bi,\bp}$ with keeping the combinatorics.
This is done in the same way as in $\Cmb_{211}$.
Hence $\mcR(\Cmb_{213})$ is connected.

\subsection{\texorpdfstring{$\Cmb_{221}$}{}}\label{subsec:221}
Let $B = C_1 + C_2 + D + M \in \mcR(\Cmb_{221})$.
As we have seen in Subsection~\ref{subsec:213}, $C_1 + C_2 + M$ can be continuously deformed to $C_{o1} + C_{o2} + L_{bi,\bp}$ with keeping the combinatorics.
Since $D$ is a member of the pencil generated by $C_1$ and $C_2$, such conic is deformed simultaneously with keeping $\Cmb_{221}$. 
Hence we infer that $B$ is continuously deformed to $C_{o1} + C_{o2} + C' + L_{bi,\bp}$, where $C'$ is a member of the pencil generated by $C_{o1}$ and $C_{o2}$.
Hence $\mcR(\Cmb_{221})$ is connected.

\subsection{\texorpdfstring{$\Cmb_{222}$}{}}\label{subsec:222}

For $\Cmb_{222}$, any element $B = C_1 + C_2 + D + M \in \mcR(\Cmb_{222})$ is determined by $C_1 + C_2 + D$.
As we have seen in \cite[Lemma~3.1]{absst},
$\mcR(\Cmb_{C_1+ C_2+ D})$
is connected and so is $\mcR(\Cmb_{222})$.

\subsection{\texorpdfstring{$\Cmb_{223}$}{}}\label{subsec:223}
This case was discussed in \cite{absst} and $\mcR(\Cmb_{223})$ has
just two connected components.

\subsection{\texorpdfstring{$\Cmb_{224}$}{}}
In \cite{tokunaga14}, we have seen there exists a Zariski pair for $\Cmb_{224}$.
Hence $\mcR(\Cmb_{224})$ has at least two connected components.
In this subsection, we will show that $\mcR(\Cmb_{224})$ has just two connected components.
We denote a member of $\mcR(\Cmb_{224})$ by $B = B_{\mcP_2} + D+ M$, where $D$ is a conic of type $D(2,4)$.
Let us start with an explicit example, which
play roles as  \lq  base points' in $\mcR(\Cmb_{224})$.

\begin{ex}
\label{ex:cmb224}
Let $\mcQ_{\btau_o}$ and $S_{\mcQ_{\btau_o}, z_o}$ ($\btau_o = (\lambda_o, \ba)$) be
 the quartic  and  the rational elliptic surface considered in Example~\ref{ex:cmb212a}.
In this case, we have
\begin{eqnarray*}
P_{12} &=& \left( -3t-2, -i \, \sqrt{10} t^{2} - 3 i \, \sqrt{10} t - 2 i \, \sqrt{10} \right), \\
P_{13} &=& \left( -t+2, -i \, \sqrt{2} t^{2} - i \, \sqrt{2} t + 2 i \, \sqrt{2} \right), \\
P_{14} &=& \left( 4,  -i \, t^{2} + 4 i \right).
\end{eqnarray*}
and
\begin{eqnarray*}
{[2]}P_{12} &=& \left( \frac{1}{10}t^2, -\frac{3}{100}i\sqrt{10}(t^2 + 20)t \right), \\
{[2]}P_{13} &=& \left( \frac{1}{2}t^2, -\frac{1}{4}i\sqrt{2}(t + 2)(t - 2)t \right), \\
{[2]}P_{14} &=& \left( t^2-\frac{9}{4}, -\frac{3}{2}i(t^2 + \frac{19}{4}) \right).
\end{eqnarray*}
Now put
\begin{eqnarray*}
D_{24,\btau_o}&:& \bar{f}_{\mcQ_{\btau_o,z_o}}(s_{[2]P_{12}}) = x-\frac{1}{10}t^2=0, \\
L_{12,\btau_o}&:& \bar{f}_{\mcQ_{\btau_o,z_o}}(s_{P_{12}}) = x+3t+2 = 0, \\
L_{13,\btau_o}&:& \bar{f}_{\mcQ_{\btau_o,z_o}}(s_{P_{13}}) = x+t-2 = 0, \\
L_{14,\btau_o}&:& \bar{f}_{\mcQ_{\btau_o,z_o}}(s_{P_{14}}) = x-4 = 0.
\end{eqnarray*}
We define $B_{1,\btau_o},B_{2,\btau_o}$ and $B_{3,\btau_o}$ to be
\[
B_{1,\btau_o} = \mcQ_{\btau_o} + D_{24,\btau_o} + L_{12,\btau_o}, \, 
B_{2,\btau_o} = \mcQ_{\btau_o} + D_{24,\btau_o} + L_{13,\btau_o}, \, 
B_{3,\btau_o} = \mcQ_{\btau_o} + D_{24,\btau_o} + L_{14,\btau_o}.
\]
Then by \cite{tokunaga10}, $(B_{1,\btau_o},B_{2,\btau_o})$ and $(B_{1,\btau_o}, B_{3,\btau_o})$ are Zariski pairs.
\end{ex}

\begin{prop}
\label{prop:cmb224}
There exists a homeomorphism $h:\PP^2 \to \PP^2$ such that $h(B_{2,\btau_o}) = B_{3,\btau_o}$. 
\end{prop}
\proof
Let $C_{1a}+C_{2a}+D_{o}$ be the one in Subsection \ref{subsec:225}.
By our argument in Subsection \ref{subsec:225}, $\mcQ_{\btau_o}+D_{24,\btau_o}$ is continuously deformed to $C_{1a}+C_{2a}+D_{o}$ keeping with $\Cmb_{\mcQ_{\btau_o}+D_{24,\btau_o}}$ such that the points $p_i \in C_1 \cap C_2$ go to $p_{j,a} \in C_{1a} \cap C_{2a}$.
Here we label $p_{j,a}$'s counterclockwisely so that $p_1$ goes to $p_{1,a}$.
Let $L_{j,a}$ be lines pass though $p_{1,a}$ and another point $p_{j,a}$ in $C_{1a} \cap C_{2a}$.
Since there is a projective transformation $\phi'$ such that $\phi'(L_{2,a})=L_{4,a}$ and $\phi'(C_{1a} + C_{2a} + D_o) = C_{1a} + C_{2a} + D_o$,  there exists a homeomorphism
$h': (\PP^2, C_{1a} + C_{2a} + D_o + L_{2,a}) \to (\PP^2, C_{1a} + C_{2a} + D_o + L_{4,a})$.

Now we show that  $L_{12,\btau_o} \rightsquigarrow L_{3,a}$. 
In fact, suppose that  $L_{12,\btau_o} \rightsquigarrow L_{2,a}$. As either 
 $L_{13,\btau_o} \rightsquigarrow L_{4,a}$ or $L_{14,\btau_o} \rightsquigarrow L_{4,a}$, this means
 that there exists a homemorphism
 from $(\PP^2, B_{1, \btau_o})$
 to $(\PP^2, B_{2, \btau_o})$ or
 $(\PP^2, B_{3, \btau_o})$, but 
 this is impossible. $L_{12,\btau_o} \rightsquigarrow L_{4,a}$ is also impossible similarly.
 Hence
 $L_{12,\btau_o} \rightsquigarrow 
 L_{3, a}$.
Thus $\{L_{13, \btau_o}, L_{14,\btau_o}\} \rightsquigarrow \{L_{2,a}, L_{4, a}\}$. Therefore
our statement follows. \endproof
 
\begin{prop}
Let $B$ be an arbitrary member in $\mcR(\Cmb_{224})$.
Then $B$ is continuously deformed to either $B_{1,\btau_o}$ or $B_{2,\btau_o}$ in Example \ref{ex:cmb224}.
In particular, $\mcR(\Cmb_{224})$ has just two connected components.
\end{prop}
\proof
After taking a suitable coordinate change, we may assume that $B$ is given as follows:

There exists $\btau \in \mcM_2$ such that
\[
B_{\btau} = \mcQ_{\btau} + D_{\btau} + L_{12,\btau}
\]
where $D_{\btau}$ is either $\tilde{f}_{\mcQ_{\btau}, z_o}(s_{[2]P_{12, \btau}})$,
$\tilde{f}_{\mcQ_{\btau}, z_o}(s_{[2]P_{13, \btau}})$ or $\tilde{f}_{\mcQ_{\btau}, z_o}(s_{[2]P_{23, \btau}})$.

\underline{Case $D_{\btau}$=$\tilde{f}_{\mcQ_{\btau}, z_o}(s_{[2]P_{12, \btau}})$.}
Consider a path $\gamma:[0,1] \to \mcM_2$ such that 
(i) $\gamma(0)=\btau$,$\gamma(1)=\btau_o$,
(ii) $B_{\gamma(s)} \in \mcR(\Cmb_{224})$ for $\forall s \in [0,1]$, and
(iii) $t_1 \rightsquigarrow -2$, $t_2 \rightsquigarrow -1$, $t_3 \rightsquigarrow 1$,
$t_4 \rightsquigarrow 2$.
Then shows that $B$ is continuously deformed to $B_{1,\btau_o}$ keeping the combinatorics.

\underline{Case $D_{\btau}$=$\tilde{f}_{\mcQ_{\btau}, z_o}(s_{[2]P_{13,\btau}})$.}
Consider a path $\gamma:[0,1] \to \mcM_2$ such that 
(i) $\gamma(0)=\btau$,$\gamma(1)=\btau_o$, 
(ii) $B_{\gamma(s)} \in \mcR(\Cmb_{224})$ for $\forall s \in [0,1]$, and
(iii) $t_1 \rightsquigarrow -2$, $t_2 \rightsquigarrow 1$, $t_3 \rightsquigarrow -1$,
$t_4 \rightsquigarrow 2$.
Note that  such that $p_2$ and $p_3$ are interchanged under
this operation.
Then $L_{12,\btau}$  
 (resp.$D_{\btau}$) is deformed to $L_{13,\btau_o}$     
 (resp.$D_{\btau_o}$).
Hence $B$ is continuously deformed to $B_{2,\btau_o}$ keeping the combinatorics.

\underline{Case $D_{\btau}$=$\tilde{f}_{\mcQ_{\btau}, z_o}(s_{[2]P_{23, \btau}})$.}
Consider a path $\gamma:[0,1] \to \mcM_2$ such that 
(i) $\gamma(0)=\btau$,$\gamma(1)=\btau_o$
(ii) $B_{\gamma(s)} \in \mcR(\Cmb_{224})$ for $\forall s \in [0,1]$,and
(iii) $t_1 \rightsquigarrow -2$, $t_2 \rightsquigarrow 2$, $t_3 \rightsquigarrow 1$, $t_4 \rightsquigarrow -1$.
Note that  such that $p_2$ and $p_4$ are interchanged under
this operation.
Then $L_{12,\btau}$ is deformed to $L_{14,\btau_o}$.
Since $[2]P_{34,\btau_o} = [2]P_{12,\btau_o}$, $D_{\btau}$ is deformed to $\tilde{f}_{\mcQ_{\btau_o}, z_o}(s_{[2]P_{34, \btau_o}})=\tilde{f}_{\mcQ_{\btau_o}, z_o}(s_{[2]P_{12, \btau_o}})=D_{\btau_o}$.
Hence $B$ is continuously deformed to $B_{3,\btau_o}$ keeping the combinatorics.

By the Proposition \ref{prop:cmb224}, $B_{2,\btau_o}$ and $B_{3,\btau_o}$ are transformed by some projective transformation each other.
Hence our statement follows. \qed

\subsection{\texorpdfstring{$\Cmb_{225}$}{}}\label{subsec:225}

Take $B = C_1 + C_2 + D + M \in \mcR(\Cmb_{225})$ arbitrary. We label four tangent points between $(C_1+C_2)\cap D$ by
$C_1\cap D = \{q_1, q_3\}$ and $C_2\cap D  = 
\{q_2, q_4\}$. $L_{q_1q_3}$
(resp. $L_{q_2q_4}$) denotes a line connecting
$q_1$ and $q_3$ (resp. $q_2$ and $q_4$). 
Then $C_1$ (resp. $C_2$) is a member of the pencil generated by
$D$ and $2L_{q_1q_2}$ (resp. $D$ and $2L_{q_2q_4}$).

Now consider a projective transformation $\phi: \PP^2 \to 
\PP^2$ such that $\phi(D) = D_o$ where $D_o$ is a conic
given by $T^2 + X^2 = Z^2$. 
Put $q_{oi} = \phi(q_i)$ ($i = 1, 2, 3, 4$). Then
$\phi(L_{q_1q_3}) = L_{q_{o1}q_{o3}}$ and
$\phi(L_{q_2q_4}) = L_{q_{o2}q_{o4}}$. Now we move
$q_{oi}$ $(i =1, 2, 3, 4)$ continuously so that
\[
q_{o1} \rightsquigarrow (1,0), \,\,
q_{o2} \rightsquigarrow (0,1), \,\,
q_{o3} \rightsquigarrow (-1,0), \,\,
q_{o4} \rightsquigarrow (0,-1).
\]
Since the two pencils of conics are also continuously
deformed along with $q_{oi}$ $(i = 1,2,3,4)$,
we infer that $C_1+ C_2 + D$ is continuously deformed
to $C_{1a} + C_{2a} + D_o$ keeping with
$\Cmb_{C_1 + C_2 +D}$, where
\[
C_{1a}: \left (\frac ta \right )^2 + x^2 = 1, \quad
C_{2a}: t^2 + \left (\frac xa \right )^2 = 1 , \quad (a \in \RR_{>1}).
\]
By Remark~\ref{rem:bitangents}, $M$ is continuously
deformed to one of $x = t \pm \sqrt{a^2 + 1},
x = -t \pm \sqrt{a^2 + 1}$. Hence
$C_{1a} + C_{2a} + D_o+ \mbox{(a bitangent)}$ is
transformed to each other by some projective transformation.
Hence $\mcR(\Cmb_{225})$ is connected.

\appendix

\section{A remark on the fundamental groups}

In this section, we study the fundamental groups of the arrangements in the Zariski pairs given in Theorem~\ref{thm:main}. We calculate a presentation of the fundamental group for each case using {\tt SageMath 10.4} \cite{sagemath} and the package {\tt Sirocco} \cite{sirocco}. Then we calculate the number of epimorphisms from the fundamental groups to $S_3$, the symmetric group of degree 3,  using {\tt GAP} \cite{GAP4}. The existence of such epimorphisms implies that the group is non-abelian, and the difference in the number of epimorphisms allows us to distinguish non-isomorphic groups. 
We use the following commands 
\begin{itemize}
    \item \verb|ProjectivePlaneCurveArrangements()|

    This command constructs projective plane curve arrangements as a {\tt SageMath} object.
    \item \verb|fundamental_group()|

    This command computes the fundamental group of the projective plane curve arrangement in terms of generators and relations. The package {\tt Sirocco} must be enabled.
    
    \item \verb|meridian()|

    This command returns the information of the meridians of the irreducible components of the arrangement in terms of the generators of the fundamental groups. The package {\tt Sirocco} must be enabled.
    
    \item \verb|GQuotients()|

    This is a {\tt GAP} command that computes epimorphisms from a group to a given finite group. The output is given in terms of the images of the generators.   
\end{itemize}
and the results are summarized in the following table:

\begin{center}
    \begin{tabular}{|c|c|c|c|}
    \hline
    Combinatorics & Arrangement & abelian/non-abelian & Num. of epi. to $S_3$\\
    \hline
    \multirow{2}{*}{$\Cmb_{123}$} & $B_{1, {\bm \tau}_o}$ & non-abelian & 5\\
    \cline{2-4}
    & $B_{2, {\bm \tau}_o}$ & non-abelian & 3\\
    \hline
    \multirow{2}{*}{$\Cmb_{124}$} & $B_{1,{\bm \tau}_o}$ & non-abelian & 7\\
    \cline{2-4}
    & $B_{2, {\bm \tau}_o} $ & non-abelian & 6\\
    \hline
    \multirow{2}{*}{$\Cmb_{212}$} & $B_{13,{\bm p}}$ & non-abelian & 7\\
    \cline{2-4}
    & $B_{12,{\bm p}}$ & non-abelian &  6\\
    \hline
    \multirow{2}{*}{$\Cmb_{223}$} &  \multirow{2}{*}{$B_1, B_2$} &  \multirow{2}{*}{free abelian of rank $3$} & \multirow{2}{*}{0} \\
    &&&\\
    \hline
    \multirow{2}{*}{$\Cmb_{224}$} & $B_{1, {\bm \tau}_o}$ & non-abelian & 7\\
    \cline{2-4}
    & $B_{2, {\bm \tau}_o}$ & non-abelian & 6\\
    \hline
    \end{tabular}
\end{center}


\begin{rem}
\begin{enumerate}
\item[(i)] The fundamental group of $\Cmb_{124}$ and $\Cmb_{224}$ were computed in \cite{asstt}.
Also the fundamental group of $\Cmb_{223}$ was calculated in \cite{absst}.
\item[(ii)] For each epimorphism to $S_3$, the orders of the images of the meridians of the irreducible components can be read off from the output of \verb|GQuotients()|. We can    
construct $S_3$-covers of 
$\PP^2$ with the corresponding branch data using the methods in \cite{Tokunaga_1994}, \cite{tokunaga10} which supports the correctness of the above calculations.
\end{enumerate}
\end{rem}







\bibliographystyle{spmpsci}
\bibliography{biblio.bib}

\begin{thebibliography}{10}
\providecommand{\url}[1]{{#1}}
\providecommand{\urlprefix}{URL }
\expandafter\ifx\csname urlstyle\endcsname\relax
  \providecommand{\doi}[1]{DOI~\discretionary{}{}{}#1}\else
  \providecommand{\doi}{DOI~\discretionary{}{}{}\begingroup
  \urlstyle{rm}\Url}\fi

\bibitem{absst}
Amram, M., Bannai, S., Shirane, T., Sinichkin, U., Tokunaga, H.o.: The
  realization space of a certain conic line arrangement of degree $7$ and a
  $\pi_1$-equivalent {Z}ariski pair (2023).
\newblock ArXiv:2307.01736 To appear in Israel J. of Math.

\bibitem{amram_garber_teicher07}
Amram, M., Garber, D., Teicher, M.: Fundamental groups of tangent conic-line
  arrangements with singularities up to order 6.
\newblock Math. Z. \textbf{256}(4), 837--870 (2007).
\newblock \doi{10.1007/s00209-007-0109-4}.
\newblock \urlprefix\url{https://doi.org/10.1007/s00209-007-0109-4}

\bibitem{asstt}
Amram, M., Shwartz, R., Sinichkin, U., Tan, S.L., Tokunaga, H.o.: {Z}ariski
  pairs of conic-line arrangements of degree $7$ and $8$ via fundamental groups
  (2023).
\newblock ArXiv:2106.0350

\bibitem{amram_teicher_uludag03}
Amram, M., Teicher, M., Uludag, A.M.: Fundamental groups of some quadric-line
  arrangements.
\newblock Topology Appl. \textbf{130}(2), 159--173 (2003).
\newblock \doi{10.1016/S0166-8641(02)00218-3}.
\newblock \urlprefix\url{https://doi.org/10.1016/S0166-8641(02)00218-3}

\bibitem{abst2023-1}
Artal~Bartolo, E., Bannai, S., Shirane, T., Tokunaga, H.: Torsion divisors of
  plane curves and {Z}ariski pairs.
\newblock St. Petersburg Math. J. \textbf{34}(5), 721--736 (2023).
\newblock \doi{10.1090/spmj/1776}.
\newblock \urlprefix\url{https://doi.org/10.1090/spmj/1776}.
\newblock Translated from Algebra i Analiz {\bf 34} (2022), no. 5

\bibitem{survey}
Artal~Bartolo, E., Cogolludo, J.I., Tokunaga, H.o.: A survey on {Z}ariski
  pairs.
\newblock In: Algebraic geometry in {E}ast {A}sia---{H}anoi 2005, \emph{Adv.
  Stud. Pure Math.}, vol.~50, pp. 1--100. Math. Soc. Japan, Tokyo (2008).
\newblock \doi{10.2969/aspm/05010001}.
\newblock \urlprefix\url{https://doi.org/10.2969/aspm/05010001}

\bibitem{bannai16}
Bannai, S.: A note on splitting curves of plane quartics and multi-sections of
  rational elliptic surfaces.
\newblock Topology Appl. \textbf{202}, 428--439 (2016).
\newblock \doi{10.1016/j.topol.2016.02.005}.
\newblock \urlprefix\url{https://doi.org/10.1016/j.topol.2016.02.005}

\bibitem{BKMT22}
Bannai, S., Kawana, N., Masuya, R., Tokunaga, H.: Trisections on certain
  rational elliptic surfaces and families of {Z}ariski pairs degenerating to
  the same conic-line arrangement.
\newblock Geom. Dedicata \textbf{216}(1), Paper No. 8, 23 (2022).
\newblock \doi{10.1007/s10711-021-00672-5}.
\newblock \urlprefix\url{https://doi.org/10.1007/s10711-021-00672-5}

\bibitem{bannai-tokunaga15}
Bannai, S., Tokunaga, H.o.: Geometry of bisections of elliptic surfaces and
  {Z}ariski {$N$}-plets for conic arrangements.
\newblock Geom. Dedicata \textbf{178}, 219--237 (2015).
\newblock \doi{10.1007/s10711-015-0054-z}.
\newblock \urlprefix\url{https://doi.org/10.1007/s10711-015-0054-z}

\bibitem{bannai-tokunaga17}
Bannai, S., Tokunaga, H.o.: Geometry of bisections of elliptic surfaces and
  {Z}ariski {$N$}-plets {II}.
\newblock Topology Appl. \textbf{231}, 10--25 (2017).
\newblock \doi{10.1016/j.topol.2017.09.003}.
\newblock \urlprefix\url{https://doi.org/10.1016/j.topol.2017.09.003}

\bibitem{dimca_pokora22}
Dimca, A., Pokora, P.: On conic-line arrangements with nodes, tacnodes, and
  ordinary triple points.
\newblock J. Algebraic Combin. \textbf{56}(2), 403--424 (2022).
\newblock \doi{10.1007/s10801-022-01116-3}.
\newblock \urlprefix\url{https://doi.org/10.1007/s10801-022-01116-3}

\bibitem{garber_friedman14}
Friedman, M., Garber, D.: On the structure of fundamental groups of conic-line
  arrangements having a cycle in their graph.
\newblock Topology Appl. \textbf{177}, 34--58 (2014).
\newblock \doi{10.1016/j.topol.2014.05.013}.
\newblock \urlprefix\url{https://doi.org/10.1016/j.topol.2014.05.013}

\bibitem{garber_friedman15}
Friedman, M., Garber, D.: On the structure of conjugation-free fundamental
  groups of conic-line arrangements.
\newblock J. Homotopy Relat. Struct. \textbf{10}(4), 685--734 (2015).
\newblock \doi{10.1007/s40062-014-0081-8}.
\newblock \urlprefix\url{https://doi.org/10.1007/s40062-014-0081-8}

\bibitem{GAP4}
The GAP Group: GAP -- Groups, Algorithms, and Programming, Version 4.13.1
  (2024).
\newblock \urlprefix\url{https://www.gap-system.org}

\bibitem{kodaira}
Kodaira, K.: On compact analytic surfaces. {II}, {III}.
\newblock Ann. of Math. (2) 77 (1963), 563--626; ibid. \textbf{78}, 1--40
  (1963).
\newblock \doi{10.2307/1970500}.
\newblock \urlprefix\url{https://doi.org/10.2307/1970500}

\bibitem{macnic}
Macinic, A.: On the deletion and addtion of a conic to a free curve (2024).
\newblock ArXiv:2408.11714

\bibitem{sirocco}
Marco-Buzunariz, M.{\'A}., Rodr{\'i}guez, M.: Sirocco: A library for certified
  polynomial root continuation.
\newblock In: G.M. Greuel, T.~Koch, P.~Paule, A.~Sommese (eds.) Mathematical
  Software -- ICMS 2016, pp. 191--197. Springer International Publishing, Cham
  (2016).
\newblock \urlprefix\url{https://doi.org/10.1007/978-3-319-42432-3_24}

\bibitem{masuya24}
Masuya, R.: Geometry of weak-bitangent lines to quartic curves and sections on
  certain rational elliptic surfaces.
\newblock Hiroshima Math. J. \textbf{54}(1), 1--31 (2024).
\newblock \doi{10.32917/h2021060}

\bibitem{miranda-BTES}
Miranda, R.: The basic theory of elliptic surfaces.
\newblock Dottorato di Ricerca in Matematica. [Doctorate in Mathematical
  Research]. ETS Editrice, Pisa (1989)

\bibitem{naruki83}
Naruki, I.: Some invariants for conics and their applications.
\newblock Publ. Res. Inst. Math. Sci. \textbf{19}(3), 1139--1151 (1983).
\newblock \doi{10.2977/prims/1195182023}.
\newblock \urlprefix\url{https://doi.org/10.2977/prims/1195182023}

\bibitem{oguiso-shioda}
Oguiso, K., Shioda, T.: The {M}ordell-{W}eil lattice of a rational elliptic
  surface.
\newblock Comment. Math. Univ. St. Paul. \textbf{40}(1), 83--99 (1991)

\bibitem{pokora-szemberg23}
Pokora, P., Szemberg, T.: Conic-line arrangements in the complex projective
  plane.
\newblock Discrete Comput. Geom. \textbf{69}(4), 1121--1138 (2023).
\newblock \doi{10.1007/s00454-022-00397-6}.
\newblock \urlprefix\url{https://doi.org/10.1007/s00454-022-00397-6}

\bibitem{schenck2009}
Schenck, H., Toh\v{a}neanu, c.S.O.: Freeness of conic-line arrangements in
  {$\mathbb P^2$}.
\newblock Comment. Math. Helv. \textbf{84}(2), 235--258 (2009).
\newblock \doi{10.4171/CMH/161}.
\newblock \urlprefix\url{https://doi.org/10.4171/CMH/161}

\bibitem{shioda90}
Shioda, T.: On the {M}ordell-{W}eil lattices.
\newblock Comment. Math. Univ. St. Paul. \textbf{39}(2), 211--240 (1990)

\bibitem{shirane18}
Shirane, T.: Connected numbers and the embedded topology of plane curves.
\newblock Canad. Math. Bull. \textbf{61}(3), 650--658 (2018).
\newblock \doi{10.4153/CMB-2017-066-5}.
\newblock \urlprefix\url{https://doi.org/10.4153/CMB-2017-066-5}

\bibitem{sagemath}
{The Sage Developers}: {S}ageMath, the {S}age {M}athematics {S}oftware {S}ystem
  ({V}ersion 10.4) (2024).
\newblock {\tt https://www.sagemath.org}

\bibitem{Tokunaga_1994}
Tokunaga, H.o.: On dihedral galois coverings.
\newblock Canadian Journal of Mathematics \textbf{46}(6), 1299^^e2^^80^^931317
  (1994).
\newblock \doi{10.4153/CJM-1994-074-4}

\bibitem{tokunaga10}
Tokunaga, H.o.: Geometry of irreducible plane quartics and their quadratic
  residue conics.
\newblock J. Singul. \textbf{2}, 170--190 (2010).
\newblock \doi{10.5427/jsing.2010.2k}.
\newblock \urlprefix\url{https://doi.org/10.5427/jsing.2010.2k}

\bibitem{tokunaga14}
Tokunaga, H.o.: Sections of elliptic surfaces and {Z}ariski pairs for
  conic-line arrangements via dihedral covers.
\newblock J. Math. Soc. Japan \textbf{66}(2), 613--640 (2014).
\newblock \doi{10.2969/jmsj/06620613}.
\newblock \urlprefix\url{https://doi.org/10.2969/jmsj/06620613}

\end{thebibliography}

\noindent Shinzo BANNAI \\
Department of Applied Mathematics, Faculty of Science, \\
Okayama University of Science, 1-1 Ridai-cho, Kita-ku,
Okayama 700-0005 JAPAN \\
{\tt bannai@ous.ac.jp}\\

\noindent  Hiro-o TOKUNAGA and Emiko YORISAKI\\
Department of Mathematical  Sciences, Graduate School of Science, \\
Tokyo Metropolitan University, 1-1 Minami-Ohsawa, Hachiohji 192-0397 JAPAN \\
{\tt tokunaga@tmu.ac.jp}
\end{document}